\documentclass[11pt,a4paper,reqno]{amsart}
\usepackage[x-mac-roman]{inputenc}
\usepackage[frenchb,english]{babel}
\usepackage[T1]{fontenc}
\usepackage{mathrsfs}
\usepackage{amsfonts,amssymb,amsmath,amsgen,amsthm}
\usepackage{color}
\usepackage{amsbsy}
\usepackage{mathrsfs}
\usepackage{bbm}
\usepackage{titletoc}
\usepackage{enumerate}
\usepackage{graphicx}
\usepackage{caption}
\usepackage{subcaption}
\usepackage[all]{xy}
\usepackage{placeins}
\theoremstyle{plain}
\newtheorem{theorem}{Theorem}[section]

\newtheorem{proposition}[theorem]{Proposition}

\newtheorem{remark}[theorem]{Remark}

\catcode`@=11
\def\section{\@startsection{section}{1}%
  \z@{1.5\linespacing\@plus\linespacing}{.5\linespacing}%
  {\normalfont\bfseries\large\centering}}
\catcode`@=12
\def\CC{{\mathbb C}}
\def\RR{{\mathbb R}}
\def\ZZ{{\mathbb Z}}
\def\calO{\mathcal O}

\def\({\left(}
\def\){\right)}
\def\<{\left\langle}
\def\>{\right\rangle}

\def\eps{\varepsilon}

\DeclareMathOperator{\RE}{Re}
\DeclareMathOperator{\IM}{Im}
\DeclareMathOperator{\DIV}{div}

\numberwithin{equation}{section}

\newcommand{\be}{\begin{equation}}
\newcommand{\ee}{\end{equation}}
\newcommand{\bea}{\begin{eqnarray}}
\newcommand{\eea}{\end{eqnarray}}
\newcommand{\bee}{\begin{eqnarray*}}
\newcommand{\eee}{\end{eqnarray*}}
\newcommand{\bse}{\begin{subequations}}
\newcommand{\ese}{\end{subequations}}

\def\eps{\varepsilon}

\def\pa{\partial}

\def\pe{{\Psi^{\eps}}}
\def\ae{{A^\eps}}
\def\se{{S^\eps}}
\def\aeu{{A^\eps_1}}
\def\aed{{A^\eps_2}}
\def\ao{{A^0}}
\def\so{{S^0}}

\def\s{{S}}
\def\a{{A}}

\begin{document}

\title[]{Uniformly accurate time-splitting methods for the semiclassical Schr\"{o}dinger equation\\
Part 1 : Construction of the schemes and simulations}
\author[P. Chartier]{Philippe Chartier}
\email{Philippe.Chartier@inria.fr}
\address{INRIA Rennes, IRMAR and ENS Rennes, IPSO Project Team, Campus de Beaulieu, F-35042 Rennes}
\author[L. Le Treust]{Lo\"{i}c Le Treust}
\email{loic.letreust@univ-rennes1.fr}
\address{IRMAR, Universit\'e de Rennes 1 and INRIA, IPSO Project}
\author[F. M\'ehats]{Florian M\'ehats}
\email{florian.mehats@univ-rennes1.fr}
\address{IRMAR, Universit\'e de Rennes 1 and INRIA, IPSO Project}

\thanks{This work was supported by the ANR-FWF Project Lodiquas ANR-11-IS01-0003 and the ANR-10-BLAN-0101 Grant (L.L.T.) and by the ANR project Moonrise ANR-14-CE23-0007-01.
}
\begin{abstract}
This article is devoted to the construction of new numerical methods for the semiclassical Schr\"odinger equation. A phase-amplitude reformulation of the equation is described where the  Planck constant $\eps$ is not a singular parameter. This  allows to build splitting schemes whose accuracy is spectral in space, of up to  fourth order in time, and independent of $\eps$ before the caustics. The second-order method additionally preserves the $L^2$-norm of the solution just as the exact flow does. In this first part of the paper, we introduce the basic splitting scheme in the nonlinear case, reveal our strategy for constructing higher-order methods, and illustrate their properties with simulations. In the second part, we shall prove a uniform convergence result for the first-order splitting scheme applied to the linear Schr\"odinger equation with a potential.
\end{abstract}

\maketitle


\section{Introduction}
In this paper, we are concerned with the numerical approximation of the solution $\pe: \RR_+\times \RR^d\to \CC$, $d\geq 1$, of the  nonlinear Schr\"odinger (NLS) equation 
\be \label{eq:GPE}
i\eps \pa_t \pe =  - \frac{\eps^2}{2} \Delta \pe + |\pe|^2\pe,
\ee
where $\eps>0$ is the so-called semiclassical parameter. The initial datum is assumed to be of the form
\be\label{initdataNLS}
	\pe(0,\cdot) = A_0(\cdot)e^{iS_0(\cdot)/\eps}\quad \mbox{ such that  } \quad \|A_0\|_{L^2(\RR^d)} = 1.
\ee
Note that the $L^2$-norm, the energy and the momentum of $\pe(0,\cdot)$, namely
\begin{eqnarray}\label{eq:invariants}
&\mbox{Mass: }&\|\pe(t,\cdot)\|_{L^2(\RR^d)}^2,\\
&\mbox{Energy: }&\int_{\RR^d}\left(\eps^2|\nabla\pe(t,x)|^2 + |\pe(t,x)|^4\right)dx,\\
&\mbox{Momentum: }&\eps\IM\int_{\RR^d}\overline{\pe(t,x)}\nabla\pe(t,x)dx,
\end{eqnarray}
are all preserved by the flow of 
\eqref{eq:GPE}, whenever $\pe(0,\cdot)\in H^1(\RR^d)$. 

Owing to its numerous occurrences in a vast number of domains of applications in physics, equation \eqref{eq:GPE} has been widely studied: existence and uniqueness results can be found for instance in \cite{cazenave2003semilinear}. In particular, it has been frequently used for the description  of Bose-Einstein condensates, as well as for the study of the propagation of laser beams (see \cite{sulem} for a detailed presentation of the physical context). In the semiclassical regime where the rescaled Planck constant $\eps$ is small, its asymptotic study allows for an appropriate description of the observables of $\pe$ through the laws of hydrodynamics. We refer to \cite{bookRemi} for a detailed presentation of the semiclassical analysis of the NLS equation and to \cite{MR2805153} for a review of both theoretical and numerical issues.

\subsection{Motivation}
Generally speaking, numerical methods for equation \eqref{eq:GPE} exhibit an error of size 
$\Delta t^p/\eps^r + \Delta x^q/\eps^s$, 
where $\Delta t$ and $\Delta x$ are the time and space steps and $p,q,r,s$ strictly positive numbers. For time-splitting methods in the linear case for instance, the error on the wave function behaves like $\Delta x/\eps+\Delta t^p/\eps$ \cite{MR1880116,MR2739463}, while in the nonlinear framework, numerical experiments  \cite{MR2047194,MR2359915,MR3047949} suggest larger values of $r$ and $s$ after the appearance of caustics\footnote{We refer to \cite{MR3047949} for a study of  local error estimates for the semiclassical NLS equation.}. Even if we content ourselves with observables in the linear case\footnote{These authors performed extensive numerical tests in both linear and nonlinear cases \cite{MR1880116,MR2047194}. 
}, the error of a splitting method of Bao, Jin, and Markowich \cite{MR1880116} grows like $\Delta x/\eps+\Delta t^p$. Now, 
achieving a fixed accuracy for varying values of $\eps$ requires to keep both ratios $\Delta t / \eps^{r/p}$ and $\Delta x/ \eps^{s/q}$ constant, and  becomes {\em prohibitively costly} when $\eps \to 0$. Our aim, in this article, is thus to develop {\em new} numerical schemes that are {\em Uniformly Accurate (UA)} w.r.t. $\eps$, i.e. whose accuracy does not deteriorate for vanishing $\eps$. In other words, schemes for which $r, s=0$. This seems highly desirable as all available methods with the exception of \cite{BCM2013}, namely finite difference methods \cite{Akrivis1993,DFP1981,KOAD1993,wu1996}, splitting methods \cite{MR1921908,MR2739463,MR3047949,MR2429878,PM1990,MR842641}, relaxation schemes \cite{MR1785963} and symplectic methods \cite{MR967679} fail to be UA. 

It is the belief of the authors that, prior to the construction of UA-schemes, it is necessary to reformulate \eqref{eq:GPE} as in \cite{BCM2013} and we now describe how this can be done. 
\subsection{Reformulation of the problem}
In the spirit of the {\em Wentzel-Kramers-Brillouin (WKB)} techniques, we decompose $\pe$ as the product of a slowly varying amplitude and a fast oscillating factor\footnote{Considering 
the WKB-ansatz \eqref{eq:WKBansatz} transforms the invariants \eqref{eq:invariants} into respectively
	\be 
\|\ae\|_{L^2(\RR^d)}^2, \quad \int_{\RR^d}\left(|\eps\nabla\ae + i\ae\nabla\se|^2+ |\ae|^4\right)dx \; \mbox{ and } \; 
\IM\int_{\RR^d}\overline{\ae}\(\eps\nabla\ae + i\ae\nabla\se\)dx.
	\ee
}
\be\label{eq:WKBansatz}
	\pe(t,\cdot) = \ae(t,\cdot)e^{i\se(t,\cdot)/\eps}.
\ee
From this point onwards, various choices are possible, depending on whether $\ae$ is complex or not\footnote{The Madelung transform \cite{madelung} relates the semiclassical limit of \eqref{eq:GPE} to hydrodynamic equations 
	\begin{equation*}
		\pe(t,\cdot) = \sqrt{\rho^\eps(t,\cdot)}e^{i\se(t,\cdot)/\eps}
	\end{equation*}
	and amounts to choosing $\ae\in \RR_+$.
	However, this formulation leads to both analytical and numerical difficulties 
	in the presence of vacuum, \textit{i.e.} whenever $\rho^\eps$ vanishes
	\cite{CDS2012,DGM2007}.
}: taking $\ae\in \CC$ leads to the following system \cite{Grenier1998} 
\bse\label{eq:SCv1}
	\begin{align}
	&\label{eq:SCv1phase}\pa_t \se + \frac{|\nabla \se|^2}{2}  + |\ae|^2= 0,\\
	&\label{eq:SCv1amplitude}\pa_t \ae +\nabla \se\cdot \nabla \ae +\frac{\ae}{2}\Delta \se
	 = \frac{i\eps\Delta \ae}{2}
	 \end{align}
\ese
with $ \se(0,\cdot) = S_0(\cdot)$ and $\ae(0,\cdot) = A_0(\cdot)$, whose analysis  leans on symmetrizable quasilinear hyperbolic systems (the first existence and uniqueness result has been obtained by Grenier \cite{Grenier1998}). Under appropriate smoothness assumptions, $(\ae, \se)  \in \CC \times \RR$ converges when  $\eps \to 0$ to the solution $(\ao, \so)$ of
\bse\label{eq:SCv2}
	\begin{align}
	&\label{eq:SCv2phase}\pa_t \so+ \frac{|\nabla \so|^2}{2}  + |\ao|^2= 0,\\
	&\label{eq:SCv2amplitude}\pa_t \ao +\nabla \so\cdot \nabla \ao +\frac{\ao}{2}\Delta \so
	 = 0.
	 \end{align}
\ese
Notice that $(\rho, v) = (|\ao|^2,\nabla \so)$ is then solution of the compressible Euler system %
\bse\label{eq:Euler}
	\begin{align}
	&\label{eq:Eulerphase}\pa_t v + v\cdot \nabla v  + \nabla\rho= 0,\\
	&\label{eq:Euleramplitude}\pa_t \rho +\DIV(\rho v) = 0.
	 \end{align}
\ese

Now, an important drawback of \eqref{eq:SCv1} stems from the formation of caustics in finite time \cite{bookRemi}:  the solution of \eqref{eq:SCv1} may indeed cease to be smooth even though $\pe$ is globally well-defined for $\eps>0$. In order to obtain global existence for $\eps>0$ (at least in the 1D-case), Besse, Carles and Méhats \cite{BCM2013} suggested an alternative formulation by introducing an asymptotically-vanishing viscosity term in the eikonal equation \eqref{eq:SCv1phase}. Therein, system \eqref{eq:SCv1} is replaced by
\bse\label{eq:SC}
	\begin{align}
	&\label{eq:SCphase}\pa_t \se + \frac{|\nabla \se|^2}{2}  + |\ae|^2= \eps^{2}\Delta \se,\\
	&\label{eq:SCamplitude}\pa_t \ae +\nabla \se\cdot \nabla \ae +\frac{\ae}{2}\Delta \se
	 = \frac{i\eps\Delta \ae}{2} -i\eps\ae\Delta\se
	 \end{align}
\ese
where $\se(0,x) = S_0(x)$, $\ae(0,x) = A_0(x)$ and where $x\in \RR^d$.
Let us emphasize that both  \eqref{eq:SCv1} and \eqref{eq:SC} are equivalent to \eqref{eq:GPE} in the following sense: as long as the solution $(\se,\ae)$ of \eqref{eq:SCv1} (resp. \eqref{eq:SC}) is smooth, the function $\pe$ defined by \eqref{eq:WKBansatz} solves \eqref{eq:GPE}. The following existence and uniqueness result given from \cite{Grenier1998,BCM2013}
is thus perfectly satisfactory for our purpose. 
\begin{theorem}[Grenier, Besse-Carles-Méhats\footnote{Besse \textit{et al.} \cite{BCM2013} also proved that in the one-dimensional case, the solution of \eqref{eq:SC} is global in time under the assumptions of Theorem \ref{thm:Grenier} and that $(\se,\ae)$ is globally uniformly bounded.
}]\label{thm:Grenier}
	Let $\eps_{\rm max}>0$. Assume that $(S_0,A_0)$ belongs to $H^{s+2}(\RR^d)\times H^s(\RR^d)$ where $s>d/2+1$. Then, there exist $T>0$, independent of $\eps \in(0,\eps_{\rm max}]$ and a unique solution
	\[
		(\se,\ae)\in C([0,T]; H^{s+2}(\RR^d)\times H^s(\RR^d))
	\]
	of system of equations \eqref{eq:SCv1}, resp. \eqref{eq:SC}. Moreover, $(\se,\ae)$ is bounded in 
	\[
		C([0,T]; H^{s+2}(\RR^d)\times H^s(\RR^d))
	\]
	 uniformly in $\eps$ and $\pe(t,\cdot) = \ae(t,\cdot)e^{i\se(t,\cdot)/\eps}$
	is the unique solution of \eqref{eq:GPE} with the initial datum \eqref{initdataNLS}.
\end{theorem}
The main advantage of the WKB reformulation \eqref{eq:SC} over \eqref{eq:GPE} is apparent: the semiclassical parameter $\eps$ does not give rise to singular perturbations\footnote{
The Cole-Hopf transformation \cite[Section $4.4.1$]{evans1998partial}
\begin{eqnarray} \label{rem:cole-hopf}
w^\eps = \exp\(-\frac{\se}{2\eps^{2}}\) - 1
\end{eqnarray}
transforms \eqref{eq:SCphase} into $\pa_tw^\eps - \frac{|\ae|^2}{2\eps^2}(w^\eps+1) = \eps^2\Delta w^\eps$ 
for which the regularizing effect of the viscosity term becomes arguably more apparent.}. 
Hence, it constitutes a good basis for the development of UA schemes (at least prior to the appearance of caustics), as witnessed by the methods introduced later in this paper.  
%
\subsection{Content of the paper}
First and only (up to our knowledge) UA schemes are based on the formulation \eqref{eq:SC} introduced in \cite{BCM2013}. Nevertheless, these schemes are still subject to CFL stability conditions and are of low order in time and space. In this paper, we consider, in lieu of finite differences as in \cite{BCM2013}, time-splitting methods, for they enjoy the following favorable features: 
\begin{enumerate}[(i)]
	\item they do not suffer from stability restrictions on the time step;
	\item they are easy to implement;
	\item they can be composed to attain high-order of convergence in time while remaining spectrally convergent in space.  
\end{enumerate}
The third point requires further explanation: standard compositions (for orders higher or equal to $3$) are inappropriate here, as they necessarily involve negative coefficients \cite{MR2134093}. At the same time, parabolic terms prevent the corresponding equations from being solved backward in time. Negative coefficients are thus forbidden. To bypass this apparent contradiction, it has become customary to resort to {\em complex} coefficients with positive real parts \cite{BCCM11-2,MR3042575,CCDV2009,MR2545819}. Now, the situation is here rendered even more involved by the presence of Schr\"odinger terms, which are incompatible with complex coefficients. The way out will consist in distributing the real and complex coefficients to the different parts of the splitting in a clever way. The strategy will be discussed thoroughly in Section \ref{sec:scheme4}.

The outline of the remaining of the paper is as follows. In Section \ref{sec:scheme12}, we will introduce first and second order (in time) splitting schemes, preserving exactly the $L^2$-norm, and in Section \ref{sec:num}, we will present extensive numerical experiments comparing our methods to the Strang splitting method studied in \cite{MR1880116,MR2047194}. The analysis of the uniform order of convergence of the splitting scheme is postponed to Part 2 of this work and, given the technicality of the proofs, will be envisaged only for order one and for the linear Schr\"odinger equation with a potential. 
%
%

%
%
\section{Second-order numerical scheme}\label{sec:scheme12}
The UA scheme that we now introduce is built upon splitting techniques (see for instance \cite{McLQ2002} for a general exposition). The first building block is the resolution of the eikonal equation. In the sequel, $h$ denotes the time step.
\subsection{The eikonal equation} 
\label{sec:phase}

In this section, we introduce two numerical schemes for equation \eqref{eq:phase2}.

\subsubsection{A semi-Lagrangian scheme for the eikonal equation}\label{sec:semilag}
	Let us rewrite equation \eqref{eq:phase2}
	\[
		\pa_t S + H(\nabla S) = 0, 
	\] 
	where $H(p) = |p|^2/2$. Thanks to the method of characteristics, we get for $x_0\in\RR^d$  and $h\geq 0$ small enough that
	\[
		S(h,x(h,x_0)) = S(0,x_0) + h\frac{|p(h,x_0)|^2}{2}
	\]
	where $(x,p)$ solves the Hamilton equation, (see \cite[Section 3.2]{evans1998partial})
	\[\begin{split}
		&\pa_h  x(h,x_0) = \pa_p H(p(h,x_0)) =  p(h,x),\\
		&\pa_h p(h,x_0) = -\pa_x H( p(h,x_0)) = 0
	\end{split}\]
	with $x(0,x_0) = x_0$ and $p(0,x_0)  = \nabla S(0,x_0)$, \textit{i.e}
	\[
		x(h,x_0) = x_0 + h\nabla S(0,x_0)  \mbox{ and } p(h,x_0) = \nabla S(0,x_0).
	\]
	Hence, defining implicitly for any $x\in \RR^d$, the function $y(h,x) = \Gamma( x,h,y(h,x))$ where
	\[
		\Gamma( x,h,y(h,x)) = x - h\nabla S(0,y(h,x)),
	\]
	we get
	\[
		S(h,x) = S(0,y(h,x)) + h\frac{|\nabla S(0,y(h,x))|^2}{2}.
	\]
	Let us define
	\[\begin{split}
		&y^1(h,x) = \Gamma( x,h,x),\\
		&y^2(h,x) = \Gamma (x,h,\Gamma(x,h,x)).
	\end{split}\]
	and for $i = 1,2$,
	\[\begin{split}
		&S^i(h,x) = S(0,x - h\nabla S(0,y^i(h,x))) + h\frac{|\nabla S(0,y^i(h,x))|^2}{2}.
	\end{split}\]
	After long but straightforward calculations, we can prove that $S^1$ and $S^2$ define methods of order $2$ and $4$  in time respectively that do not suffer from CFL restrictions.
	\begin{remark}
		The authors think that $y^k(h,x) = \Gamma( x,h,y^{k-1}(h,x))$ gives rise to a method of order $2k$ in time.
	\end{remark}
	In the case where the space is also discretized, the functions $S(0,\cdot)$ and $\nabla S(0,\cdot)$ are interpolated thanks to the discrete Fourier transform. 
	
	Most of the numerical schemes introduced for equation \eqref{eq:phase2} (ENO and WENO type or semi-lagrangian methods) are designed to deal with low regularity solutions. In general, their precisions is of low order in time and space and the procedure to get high order schemes is very involved. The finite difference methods even suffer from CFL restrictions. 
	
	The methods we introduced are of spectral precision in space and can be of any order in time. This can be achieved in our case since we deal with regular solutions.

\subsubsection{A splitting scheme for the eikonal equation}\label{sec:splieik}

Let us remind that the Cole-Hopf transformation 
\[
	w^\eps = \exp\(-\frac{\se}{2\eps^{2}}\) - 1
\]
ensures that solving 
\be\label{eq:phase0}
	\pa_t\se + \frac{|\nabla \se|^2}{2} = \eps^{2}\Delta\se,\quad  \se(0,\cdot) = S_0(\cdot),
\ee
is equivalent to solving the heat equation
\be\label{eq:phase1}
	\pa_tw^\eps =  \eps^{2}\Delta w^\eps,\quad w^\eps(0,\cdot) = \exp\(-\frac{S_0(\cdot)}{2\eps^{2}}\) - 1,
\ee
(see Remark \ref{rem:cole-hopf}).

Due to the limitations of the floating-point representation, it is completely inaccurate to solve \eqref{eq:phase0} using \eqref{eq:phase1} for small values of $\eps$. To overcome this difficulty, we split the nonviscous eikonal equation \eqref{eq:phase2} into two subequations, at each time step; the parabolic term $\eps^2\Delta \s$ will be dealt with separately. The key idea is to allow $\s$ to be a complex-valued function despite the fact that the solutions of equations \eqref{eq:SCv1phase}, \eqref{eq:SCphase}, \eqref{eq:phase0} and \eqref{eq:phase2} take their values in $\RR$. 
		\smallskip
	
		\noindent
		{\em First flow: }\label{sec:phaseflow1} let us define $\phi^1_h$ as the exact flow at time $h\in \RR$ of equation
		\be\label{eq:phaseflot1}
		\pa_t S + \frac{\nabla S\cdot \nabla S}{2} -i\Delta S = 0
		\ee
		where $\nabla S\cdot \nabla S = \sum_{k = 1}^d (\pa_kS)^2$. 
		This equation can be solved thanks to the following modified Cole-Hopf transformation
		\[
			w = \exp\(\frac{iS}{2}\) - 1
		\]
		leading to
		\be\label{eq:phaseflot1Schro}
			i\pa_t w= -\Delta w,\quad w(0,\cdot) = \exp\(\frac{iS(0,\cdot)}{2}\) - 1,
		\ee
		solved in the Fourier space and
		\[
			S(h,\cdot) - S(0,\cdot) = -2i\log\(1 + \frac{w(h,\cdot)-w(0,\cdot)}{w(0,\cdot)+1}\).
		\]
		\begin{remark}\label{rem:condeiklog}
			This formula is well-defined whenever 
			\[
				\left\|\frac{w(h,\cdot)-w(0,\cdot)}{w(0,\cdot)+1}\right\|_{L^\infty}<1,
			\]
			a condition ensured as soon as $h$ is small enough.
		\end{remark}
			\smallskip
	
	\noindent
		{\em Second flow: }\label{sec:phaseflow2} $\phi^2_h$ is the exact flow at time $h\in \RR$ of the free Schr\"odinger equation
		\be\label{eq:phaseflot2}
			\pa_t S+i\Delta S = 0.
		\ee
			\smallskip
			
		\noindent We are now able to define time-splitting schemes for equation \eqref{eq:phase2}.  In particular, the \textit{Lie-Trotter splitting formula},
		\be\label{eq:schemephase1}
			\RE \phi^1_h\circ\phi^2_h
		\ee
		gives us an approximation of first-order in time of the solution of \eqref{eq:phase2} (which is real-valued)  while the \textit{Strang splitting formula} 
		\be\label{eq:schemephase2}
			\RE \phi^1_{h/2}\circ\phi^2_h\circ \phi^1_{h/2}
		\ee
		provides a second-order method.
	\subsection{Numerical schemes for system of equations \eqref{eq:SC}}
	We are now  in position to introduce our numerical schemes. To this aim, we split system \eqref{eq:SC}
	\begin{align*}
		&\pa_t \se + \frac{|\nabla \se|^2}{2}  + |\ae|^2 = \eps^{2}\Delta \se,\\
		&\pa_t \ae +\nabla \se\cdot \nabla \ae +\frac{\ae}{2}\Delta \se
		 = \frac{i\eps\Delta \ae}{2} -i\eps\ae\Delta\se,
	\end{align*}
	into four subsystems.
		\smallskip
	
	\noindent
	{\em First flow: }  Let us define $\varphi^1_h$ as the approximate flow at time $h\in \RR$ of the system of equations: 
	\bse\label{eq:SCflot1}
		\begin{align}
		&\label{eq:SCflot1phase}\pa_t S + \frac{|\nabla S|^2}{2}  = 0,\\
		&\label{eq:SCflot1amplitude}\pa_t A +\nabla S\cdot \nabla A +\frac{A}{2}\Delta S
	 	= \frac{i\Delta A}{2}.
		\end{align}
	\ese
	The eikonal equation \eqref{eq:SCflot1phase} is solved according to Sec. \ref{sec:phase}. Equation \eqref{eq:SCflot1amplitude} is dealt with by noticing that $w = A \exp\(iS\)$ satisfies the free Schr\"odinger equation
	\[
		i\pa_t w = -\frac{1}{2}\Delta w.
	\]
		\smallskip
	
	\noindent
	{\em Second flow: } $\varphi^2_h$ is the exact flow at time $h\in\RR$ of 
	\bse\label{eq:SCflot2}
		\begin{align}
		&\label{eq:SCflot2phase}\pa_t S = 0,\\
		&\label{eq:SCflot2amplitude}\pa_t A 
	 	= \frac{i\(\eps-1\)\Delta A}{2}
		\end{align}
	\ese
	solved in the Fourier space.
		\smallskip
	
	\noindent
	{\em Third flow: }   $\varphi^3_h$ is the exact flow at time $h\in\RR$ of 
	\bse\label{eq:SCflot3}
		\begin{align}
		&\label{eq:SCflot3phase}\pa_t S  = -|A|^2, \\
		&\label{eq:SCflot3amplitude}\pa_t A  =0.
		\end{align}
	\ese
	\smallskip
	
	\noindent
	{\em Fourth flow: }   $\varphi^4_h$ is the exact flow at time $h\in\RR_+$ of 
	\bse\label{eq:SCflot4}
		\begin{align}
		&\label{eq:SCflot4phase}\pa_t S  = \eps^2\Delta S, \\
		&\label{eq:SCflot4amplitude}\pa_t A  =-i\eps A\Delta S.
		\end{align}
	\ese
	Equation \eqref{eq:SCflot4phase} is solved in the Fourier space and the solution of \eqref{eq:SCflot4amplitude} is
	\begin{align*}
		A(h,\cdot) = \exp\(-i\eps^{-1}(S(h,\cdot)-S(0,\cdot))\)A(0,\cdot).
	\end{align*}
	Remark that $\varphi^4_h$ gathers the terms of \eqref{eq:SC} which are not in \eqref{eq:SCv1} and can thus be viewed as a regularizing flow. 

		\begin{remark}\label{rem:defprop}
		Let us stress that $\varphi^4_h$ is not defined for $h$ such that $\RE h<0$. As a matter of fact, the propagator $e^{z\Delta}$ is well-defined, in the distributional sense, if and only if $\RE(z)\geq0$.
		\end{remark}
	\smallskip
	We consider now the following methods for  \eqref{eq:SC}
	\be\label{scheme1}
		 \varphi^1_h\circ\varphi^2_h\circ\varphi^3_h\circ \varphi^4_h
	\ee
	 and 
	 \be\label{scheme2}
	 	 \varphi^1_{h/2}\circ\varphi^2_{h/2}\circ\varphi^3_{h/2}\circ \varphi^4_{h}\circ \varphi^3_{h/2}\circ\varphi^2_{h/2}\circ\varphi^1_{h/2}.
	 \ee

	\begin{remark}
		It is worth mentioning that both schemes preserve exactly the $L^2$ norm of $\a$ since all $\varphi^1_h$, $\varphi^2_h$, $\varphi^3_h$ and $\varphi^4_h$ do so.
	\end{remark}
	\section{Fourth-order numerical scheme}\label{sec:scheme4}
		The splitting of \eqref{eq:SC} into the four flows \eqref{eq:SCflot1}, \eqref{eq:SCflot2},  \eqref{eq:SCflot3}, \eqref{eq:SCflot4} proposed in the previous section is incompatible with splitting methods of order higher than $2$ with real-valued coefficients. Indeed, it is known that such methods involve at least one negative time step for each part of the splitting (see for instance \cite{MR2134093}). Therefore, we cannot built such a scheme for \eqref{eq:SC}  because of its time irreversibility.
				
		To circumvent this difficulty, it is possible to use splitting methods with complex coefficients \cite{CCDV2009,BCCM11-2,MR3042575,MR2545819}. Let us point out the main restrictions on the coefficients in order for the methods to be well-defined. For obvious consistency reasons, if a flow $\varphi_h$ is used with complex coefficients, then both $\varphi_{\alpha h}$ and $\varphi_{\beta h}$ with $\RE \alpha>0$ and $\IM \beta<0$ will appear. Hence, the flows $\varphi^1_h$ of \eqref{eq:SCflot1} and $\varphi^2_h$ of \eqref{eq:SCflot2} containing Schr\"odinger type terms have to be integrated with only real-valued coefficients, otherwise some parts of the splitting would be ill-posed. Moreover, $\varphi^4_h$ of \eqref{eq:SCflot4} contains parabolic terms and it should be used with coefficients with nonnegative real part.  In this section, we introduce a four-flow complex splitting method taking into account all these constraints. 
		
		 The main remaining problem originates from the non-analytic character of the nonlinearity appearing in the flow $\varphi^3_h$ of \eqref{eq:SCflot3}. To overcome it, we split the real and imaginary parts of $\ae$  as in \cite{BCCM11-2}.
		
		Although we content ourselves in the sequel with a fourth-order scheme, let us emphasize that the strategy adopted here is amenable to higher orders.

	\subsection{The new splitting scheme}
	We commence from the original system \eqref{eq:SC}
	\begin{align*}
		&\pa_t \se + \frac{|\nabla \se|^2}{2}  + |\ae|^2 = \eps^{2}\Delta \se,\\
		&\pa_t \ae +\nabla \se\cdot \nabla \ae +\frac{\ae}{2}\Delta \se
		 = \frac{i\eps\Delta \ae}{2} -i\eps\ae\Delta\se
	\end{align*}
	and rewrite it (following the steps exposed in \cite{BCCM11-2}) in term of the unknowns $\aeu = \RE \ae$, $\aed = \IM \ae$ and $\se$:
		\bse\label{eq:SCcomplex}
		\begin{align}
		&\label{eq:SCcomplexphase}\pa_t \se + \frac{|\nabla \se|^2}{2}  + (\aeu)^2 + (\aed)^2 = \eps^{2}\Delta \se,\\
		&\label{eq:SCcomplexamplitude} \(\pa_t  +\nabla \se\cdot \nabla +\frac{\Delta \se}{2}\)\(\begin{array}{c}\aeu\\\aed\end{array}\) = -i\eps\sigma_2\(\frac{\Delta}{2}-\Delta \se\)\(\begin{array}{c}\aeu\\\aed\end{array}\)
		\end{align}
	\ese
	in order for the nonlinearity $(\aeu)^2 + (\aed)^2$ to be an analytic function of $\aeu$ and $\aed$.
	The matrix $\sigma_2$ is here the second Pauli matrix
	\[
		\sigma_2 = \(\begin{array}{cc}0&-i\\i&0\end{array}\),
	\]
	so that $P\sigma_2 P = \sigma_3$ with
	\be\label{eq:matP}
 \sigma_3 = \(\begin{array}{cc}1&0\\0&-1\end{array}\)\quad \mbox{and}\quad		P = \frac{1}{\sqrt{2}}\(\begin{array}{cc}1 & -i \\i&-1\end{array}\) =  \frac{1}{\sqrt{2}}\(\sigma_2+\sigma_3\).
	\ee
	Let 
	\[
		V^\eps =\(\begin{array}{c}v^\eps_1\\v^\eps_2\end{array}\) = P\(\begin{array}{c}\aeu\\\aed\end{array}\).
	\]
	System \eqref{eq:SCcomplex} becomes
	\bse\label{eq:SCcomplexv2}
		\begin{align}
		&\label{eq:SCcomplexphasev2}\pa_t \se + \frac{|\nabla \se|^2}{2}  -2iv^\eps_1v^\eps_2 = \eps^{2}\Delta \se,\\
		&\label{eq:SCcomplexamplitudev2} \(\pa_t  +\nabla \se\cdot \nabla +\frac{\Delta \se}{2}\)V^\eps = -i\eps\sigma_3\(\frac{\Delta}{2}-\Delta \se\)V^\eps.
		\end{align}
	\ese
	We are now in position to define a four-flow splitting which is compatible with complex coefficients.

	\smallskip
	
	\noindent
	{\em First flow: }  Let us define $\widetilde \varphi^1_h$ as the approximate flow at time $h\in \RR$ of: 
	\bse\label{eq:SCflot1v2}
		\begin{align}
		&\label{eq:SCflot1phasev2}\pa_t S + \frac{\nabla S\cdot \nabla S}{2}  = 0,\\
		&\label{eq:SCflot1amplitudev2}\pa_t V +\nabla S\cdot \nabla V +\frac{\Delta S}{2}V
	 	= \frac{i\Delta V}{2}.
		\end{align}
	\ese
	To solve  \eqref{eq:SCflot1phasev2}, we use either the semi-lagrangian method of order $4$ of Section \ref{sec:semilag} or the following fourth-order time-splitting from \cite{Yoshida}
	\be\label{eq:yoshida}
		\phi^1_{\alpha_1h}\circ \phi^2_{\alpha_2h}\circ \phi^1_{\alpha_3h}\circ \phi^2_{\alpha_4h}\circ \phi^1_{\alpha_5h}\circ \phi^2_{\alpha_6h}\circ \phi^1_{\alpha_7h},
	\ee
	with coefficients $\alpha_1,\dots \alpha_7$ defined by
	\be\label{eq:coefalpha}
	\begin{split}
		&\alpha_1 = \alpha_7 = \frac{1}{2(2-2^{1/3})},\qquad \alpha_3 = \alpha_5 = 0.5 - \alpha_1,\\
		 &\alpha_2 = \alpha_6= \frac{1}{2-2^{1/3}},\qquad \alpha_4 = 1 - 2\alpha_2
	\end{split}
	\ee
	 and where the numerical flows $\phi^1_h$ and $\phi^2_h$ are those introduced in Sec. \ref{sec:splieik}. To solve \eqref{eq:SCflot1amplitudev2}, we proceed as for \eqref{eq:SCflot1amplitude}.
	\smallskip
	
	\noindent
	{\em Second flow: } $\widetilde \varphi^2_h$ is the exact flow at time $h\in\RR$ of 
	\bse\label{eq:SCflot2v2}
		\begin{align}
		&\label{eq:SCflot2phasev2}\pa_t S = 0,\\
		&\label{eq:SCflot2amplitudev2}\pa_t V 
	 	= -\frac{i\(\eps\sigma_3+\sigma_0\)\Delta V}{2}
		\end{align}
	\ese
	solved in the Fourier space. Here $\sigma_0$ denotes the identity matrix.
	\smallskip
	
	\noindent
	{\em Third flow: }   $\widetilde \varphi^3_h$ is the exact flow at time $h\in\RR$ of 
	\bse\label{eq:SCflot3v2}
		\begin{align}
		&\label{eq:SCflot3phasev2}\pa_t S  = 2iv_1v_2, \\
		&\label{eq:SCflot3amplitudev2}\pa_t V  =0.
		\end{align}
	\ese
	\smallskip
	
	\noindent
	{\em Fourth flow: }   $\widetilde \varphi^4_h$ is the exact flow at time $h\in\{z\in \CC,\,\RE z\geq 0\}$ of 
	\bse\label{eq:SCflot4v2}
		\begin{align}
		&\label{eq:SCflot4phasev2}\pa_t S  = \eps^2\Delta S, \\
		&\label{eq:SCflot4amplitudev2}\pa_t V  =i\eps\sigma_3 V\Delta S.
		\end{align}
	\ese
	We will see below that in the complex splitting method that we use,  the coefficients related to $\widetilde \varphi^4_h$ are complex so that $\s$, $v_1$ and $v_2$ are complex-valued functions.
	\subsection{Splitting scheme of fourth-order for system \eqref{eq:SC}}
	We define below a complex splitting  method whose coefficients related to $\widetilde \varphi^4_h$ have positive real part whereas those associated with  $\widetilde \varphi^1_h$ and $\widetilde \varphi^2_h$ are real-valued. 
	
	The simplest way to get a fourth-order time-splitting scheme for four flows is to compose several times a fourth-order time-splitting for two flows:  
 using the same time-splitting scheme as in \eqref{eq:yoshida}, we define the following fourth-order schemes (for $h\in \RR$):
	\[
		\widetilde \varphi^{12}_h = \widetilde \varphi^2_{\alpha_1h}\circ\widetilde \varphi^1_{\alpha_2h}\circ\widetilde \varphi^2_{\alpha_3h}\circ\widetilde \varphi^1_{\alpha_4h}\circ\widetilde \varphi^2_{\alpha_5h}\circ\widetilde \varphi^1_{\alpha_6h}\circ\widetilde \varphi^2_{\alpha_7h}
	\]
	for the system of equations
		\begin{align*}
		&\pa_t \s + \frac{\nabla \s\cdot \nabla \s}{2}   = 0,\\
		& \pa_tV  +\nabla \s\cdot \nabla V +\frac{\Delta \s}{2}V = -i\eps\sigma_3\frac{\Delta V}{2}.
		\end{align*}
	and
	\[
		\widetilde \varphi^{123}_h = \widetilde \varphi^{12}_{\alpha_1h}\circ\widetilde \varphi^3_{\alpha_2h}\circ\widetilde \varphi^{12}_{\alpha_3h}\circ\widetilde \varphi^3_{\alpha_4h}\circ\widetilde \varphi^{12}_{\alpha_5h}\circ\widetilde \varphi^3_{\alpha_6h}\circ\widetilde \varphi^{12}_{\alpha_7h},
	\]
	for the system of equations
	\begin{align*}
		&\pa_t S + \frac{\nabla S\cdot \nabla S}{2} -2iv_1v_2 = 0,\\
		&\pa_t V +\nabla S\cdot \nabla V +\frac{\Delta S}{2}V
	 	= -i\eps\sigma_3\frac{\Delta V}{2}.
	\end{align*}	
	The coefficients $\alpha_1,\dots,\alpha_7$ are defined by \eqref{eq:coefalpha}.
	Since $\widetilde \varphi^4_{h}$ is not reversible, we cannot use the scheme \eqref{eq:yoshida} anymore to define our four-flow method, given that $\alpha_3$, $\alpha_4$ and $\alpha_5$ are negative. To avoid this problem, we use a complex splitting method of Blanes {\em et al.} \cite{BCCM11-2}:
	\be\label{scheme4}
		\widetilde \varphi^{1234}_h =  \widetilde P\(\widetilde \varphi^4_{\beta_1h}\circ\widetilde \varphi^{123}_{\beta_2h}\circ\widetilde \varphi^4_{\beta_3h}\circ\widetilde \varphi^{123}_{\beta_4h}\circ\widetilde \varphi^4_{\beta_5h}\circ\widetilde \varphi^{123}_{\beta_6h}\circ\widetilde \varphi^4_{\beta_7h}\circ\widetilde \varphi^{123}_{\beta_8h}\circ\widetilde \varphi^4_{\beta_9h}\)\widetilde P
	\ee
	where 
	\[
		\widetilde P = \(\begin{array}{cc}1 & 0\\0& P\end{array}\),
	\]
	$P$ is defined in \eqref{eq:matP} and
	\begin{align*}
		&\beta_1 = \beta_9 = 0.060078275263542357774 - 0.060314841253378523039i,\\
    		&\beta_2  = \beta_8 = 0.18596881959910913140,\\
   		&\beta_3 = \beta_7 = 0.27021183913361078161 + 0.15290393229116195895i,\\
    		&\beta_4 = \beta_6 = 0.5 - \beta_2 = 0.31403118040089086860,\\
    		&\beta_5 = 1-2\beta_1-2\beta_3 =  0.33941977120569372122 - 0.18517818207556687181i.
	\end{align*}
	Observe that all the coefficients $\beta_1,\beta_3,\beta_5, \beta_7$ and $\beta_9$ for the irreversible flow $\widetilde \varphi^4_{h}$ have a positive real part and all the coefficients $\beta_2,\beta_4,\beta_6$ and $\beta_8$ for the flow $\widetilde \varphi^{123}_{h}$ containing all the Schr\"odinger terms are real-valued.

	Following \cite{BCCM11-2}, we state without proof the following proposition.
	\begin{proposition}\label{prop:conv4}
	 	Assume that $(A_0,S_0)\in H^s(\RR^d)\times H^{s+2}(\RR^d)$ for a large enough $s$.  Then, the following error bound holds true 
		\[
			\begin{split}
				&\|\widetilde \varphi^{1234}_h(S_0,\RE A_0, \IM A_0 )-(\se(h,\cdot),\RE \ae(h,\cdot),  \IM \ae(h,\cdot))\|_{L^2}\leq Ch^5\\
			\end{split}
		\]
		where $C$ does not depend on $\eps\in[0,\eps_{\rm max}]$ and $(\se,\ae)$ is the solution of system of equations \eqref{eq:SC}.
	 \end{proposition}
	 Let us remark that since $\se$ takes complex values, the $L^2$ norm of $\ae$ is not exactly conserved by the flows $\widetilde \varphi^{1}_h$  and $\widetilde \varphi^{4}_h$. We stress that the function $\se$ is not  projected on the set of real-valued functions after each flow, as it was done in Section \ref{sec:splieik}, since it would reduce the order of convergence.
	\section{Numerical experiments}\label{sec:num}
	In this part, we illustrate the behavior of the schemes \eqref{scheme2} and \eqref{scheme4} introduced in Sections \ref{sec:scheme12} and \ref{sec:scheme4} and compare their properties to those of the Strang splitting method \cite{MR2047194} for \eqref{eq:GPE}. As mentioned in the introduction, quadratic observables have some peculiarities for this problem. For this reason, the convergence properties of the different schemes will be illustrated separately, on the one hand for the functions $S^\eps$, $A^\eps$ (resp. $\Psi^\eps$ for the Strang splitting scheme) and, on the other hand, for the density $\rho^\eps=|A^\eps|^2$ (resp.  $\rho^\eps=|\Psi^\eps|^2$). We restrict ourselves to the one-dimensional periodic setting in which the equations studied remain unchanged. 
	
	We consider the following initial data:
	\be\label{eq:initialdata}
	\begin{split}
		&A_0(x) = \sin(x),\quad S_0(x) ={\sin(x)}/{2},\\
	        &\pe(0,\cdot) = A_0(\cdot)e^{iS_0(\cdot)/\eps},
	\end{split}
	\ee
	where $x\in\mathbb{T} = \RR/2\pi\ZZ$, for which caustics appear numerically at time $T_c = 0.5$. In our simulations, the semiclassical parameter $\eps$ varies from $1$ to $2^{-12}$.
	
	The numerical solutions $(S^\eps, A^\eps)$, resp. $\Psi^\eps$, are compared to corresponding reference solutions $(S^\eps_{ref}, A^\eps_{ref})$, resp. $\Psi^\eps_{ref}$, which, in the absence of analytical solutions, are respectively obtained thanks to our fourth order splitting method \eqref{scheme4} and thanks to a splitting scheme of order $4$ for \eqref{eq:GPE} (see \cite{Yoshida,MR2047194}), with very small time and space steps. More precisely, to compute $(S^\eps_{ref}, A^\eps_{ref})$, we have taken $N_x = 2^{8}$ and $h = 2^{-13}T_f$, and to compute $\Psi^\eps_{ref}$, in order to fit with the constraints on the time step and on the space step
	\[
		h \ll \eps \mbox{ and } \Delta x \ll \eps,
	\]
 the space interval $[0,2\pi]$ is discretized with $N_x = 2^{15}$ points and the time step is $h = 2^{-18}T_f$. 	 

	The various errors that are represented on the figures below are defined as follows:
	\[
		\begin{split}
			&err_{\rho^\eps}(T) = \frac{\|\rho_{ref}^\eps(T) - \rho^{\eps}(T)\|_{L^1}}{\|\rho_{ref}^\eps(T)\|_{L^1}},\quad err_{\pe}(T) = \frac{\|\Psi^\eps_{ref}(T)-\pe(T)\|_{L^2}}{\|\psi^\eps_{ref}(T)\|_{L^2}},
		\end{split}
	\]
	and
	\[
		\begin{split}	
			&err_{(\se,\ae)}(T) = \(\frac{\|S^\eps_{ref}(T)-\se(T)\|_{L^2}^2+\|A^\eps_{ref}(T)-\ae(T)\|_{L^2}^2}{\|S^\eps_{ref}(T)\|_{L^2}^2+\|A^\eps_{ref}(T)\|_{L^2}^2}\)^{1/2},
		\end{split}
	\]
	where 
	\[
		\|u\|_{L^1} = \Delta x \sum_{k = 0}^{N_x-1}|u_k|, \quad\|u\|_{L^2} = \sqrt{\Delta x \sum_{k = 0}^{N_x-1}|u_k|^2},
	\] 
	and $\rho_{ref}^\eps(T) = |\Psi^{\eps}_{ref}(T)|^2$. As far as the Strang splitting scheme is concerned, $\rho^{\eps}(T) = |\Psi^{\eps}(T)|^2$ whereas $\rho^{\eps}(T) = |\ae(T)|^2$ for our methods. 

\begin{figure}[p]
	\centering
	\hspace{-4cm}
	\begin{subfigure}[t]{0.65\textwidth}
		\centering
		\includegraphics[height=5.cm,width=\textwidth]{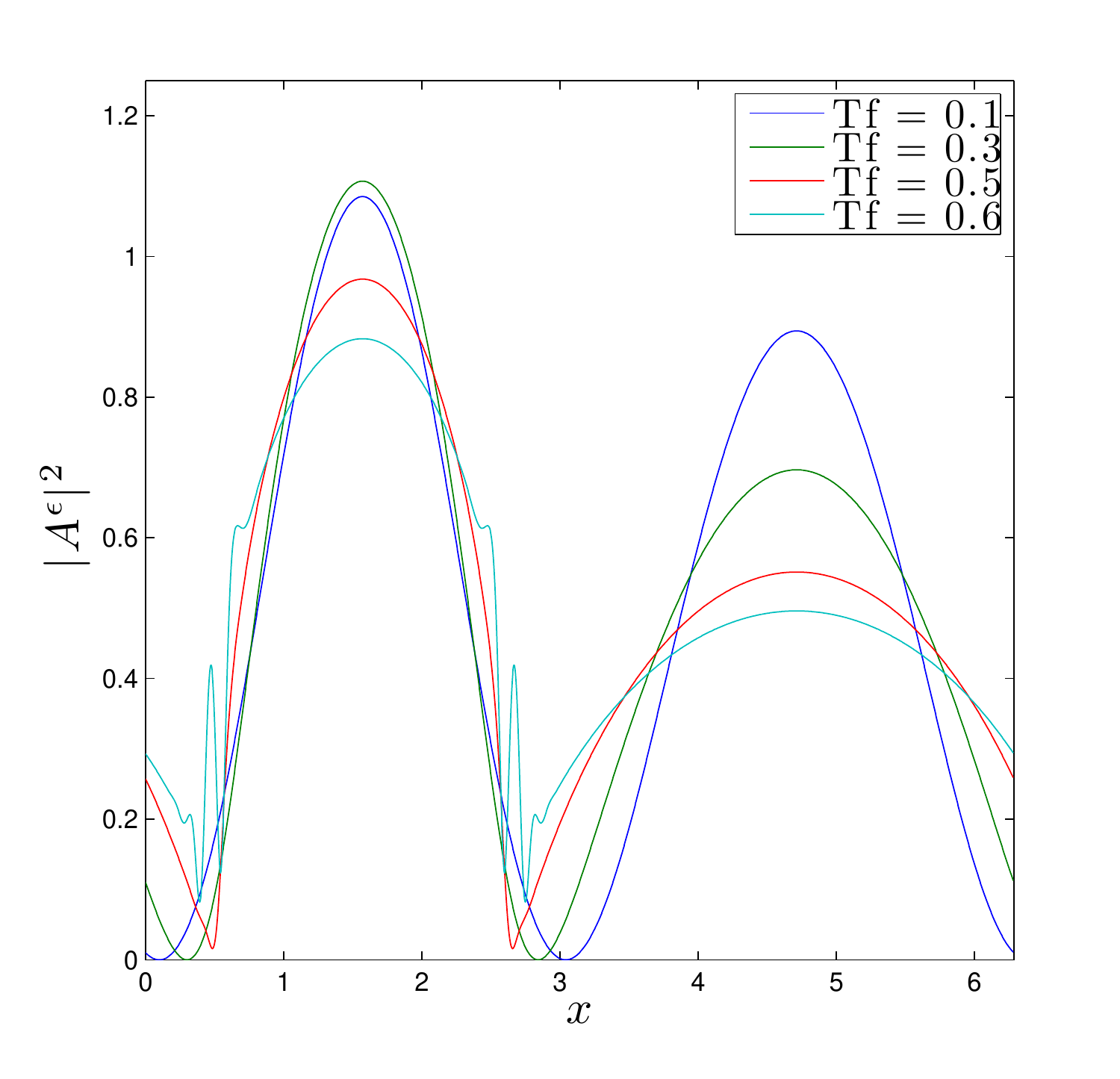}
		\caption{\label{fig:evolDensity} Evolution of the density $|A^\eps|^2$ for $\eps = 2^{-5}$.} 	
	\end{subfigure}
	\hspace{-1cm}
	\begin{subfigure}[t]{0.65\textwidth}
		\centering
		\includegraphics[height=5.cm,width=\textwidth]{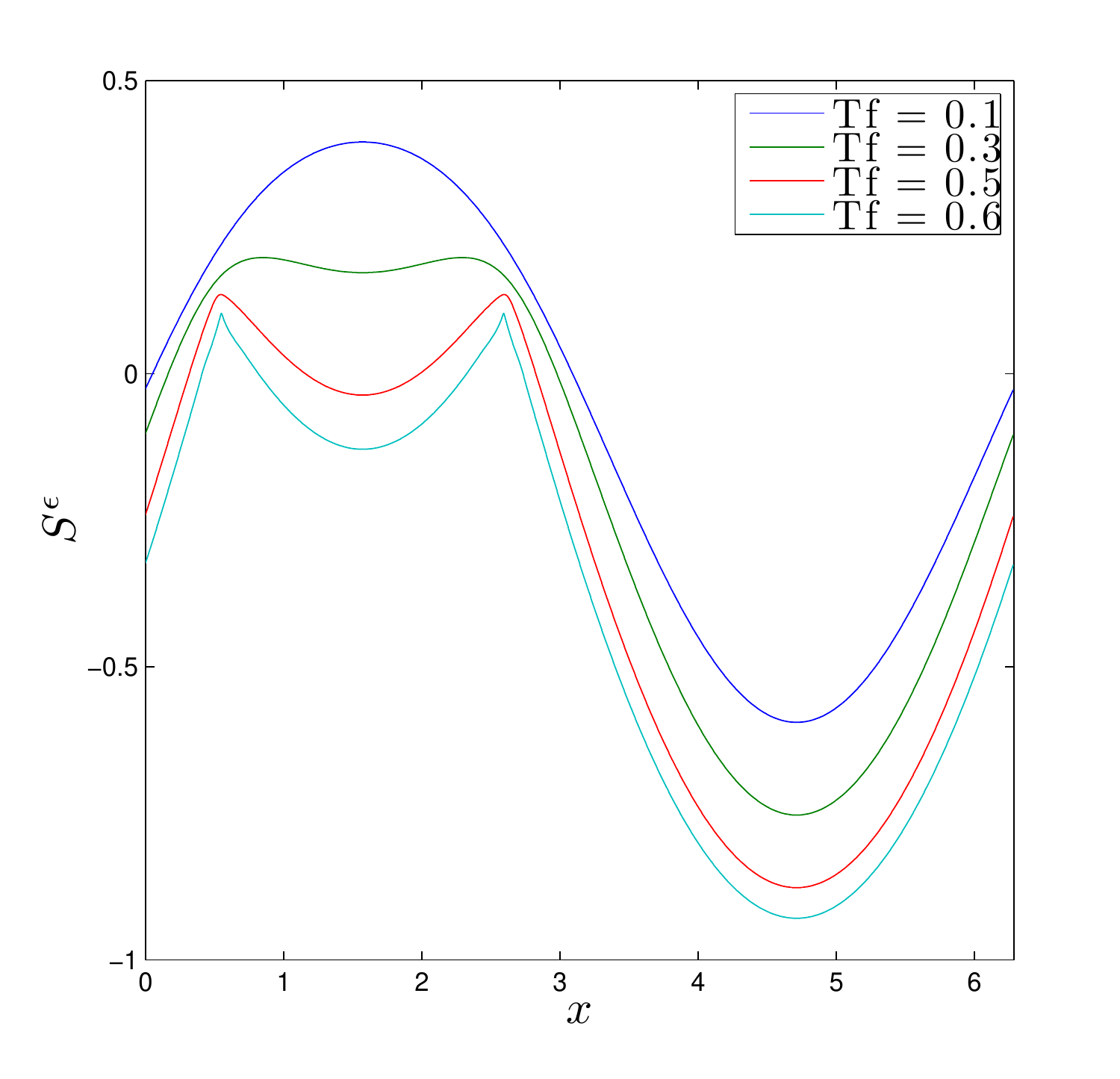}
		\caption{\label{fig:evolPhase} Evolution of the phase $S^\eps$ for $\eps = 2^{-5}$.} 
	\end{subfigure}
	\hspace{-4cm}
	\caption{Evolution of the density and of the phase}
	\end{figure}

	We first study qualitatively the dynamics, in order to guess what is the time of appearance of the caustics.
	Figures \ref{fig:evolDensity} and \ref{fig:evolPhase} represent the density $|A^\eps|^2$ and the phase $S^\eps$ at times $T_f = 0.1$, $0.3$, $0.5$, $0.6$  for $\eps = 2^{-5}$. The caustics appear around $t=0.5$. At time $t=0.6$, oscillations at other scales than those of the phase can be observed in $|A^\eps|^2$ whereas $S^\eps$ ceases to be smooth. These figures are obtained by using our scheme \eqref{scheme2} with $N_x = 2^{11}$ and $N_t = T_f/h = 2^{13}$.
	
\begin{figure}[p]
	\centering
	\hspace{-4cm}
	\begin{subfigure}[t]{0.6\textwidth}
		\centering
		\includegraphics[height=5.cm,width=\textwidth]{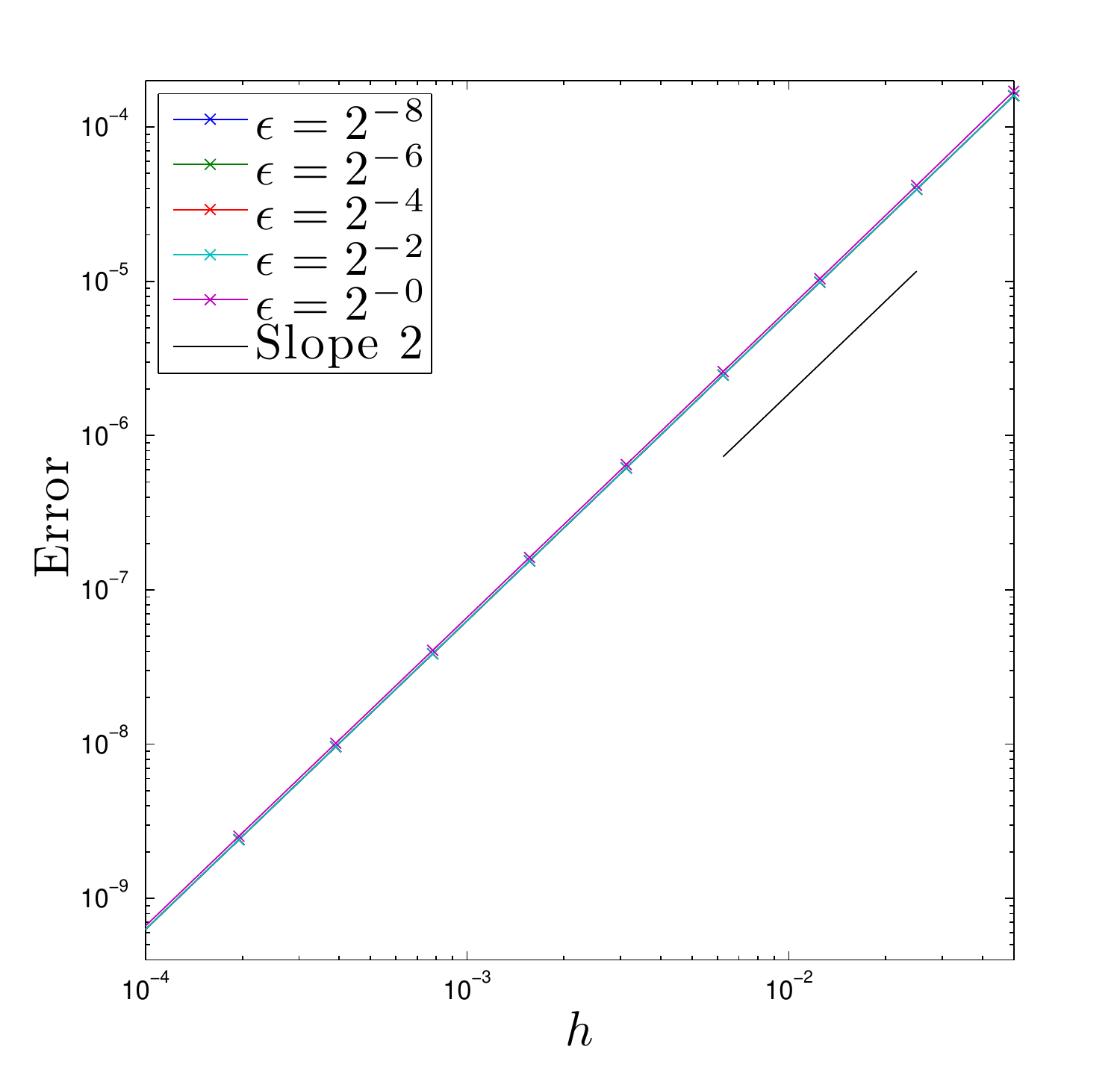}
		\caption{\label{fig:Tf1cvgtpsord2v1NLS} $err_{\rho^\eps} (T_f = 0.1)$ w.r.t $h$, $N_x = 2^{9}$}
	\end{subfigure}
	\hspace{-1cm}
	\begin{subfigure}[t]{0.7\textwidth}
		\centering
		\includegraphics[height=5.cm,width=\textwidth]{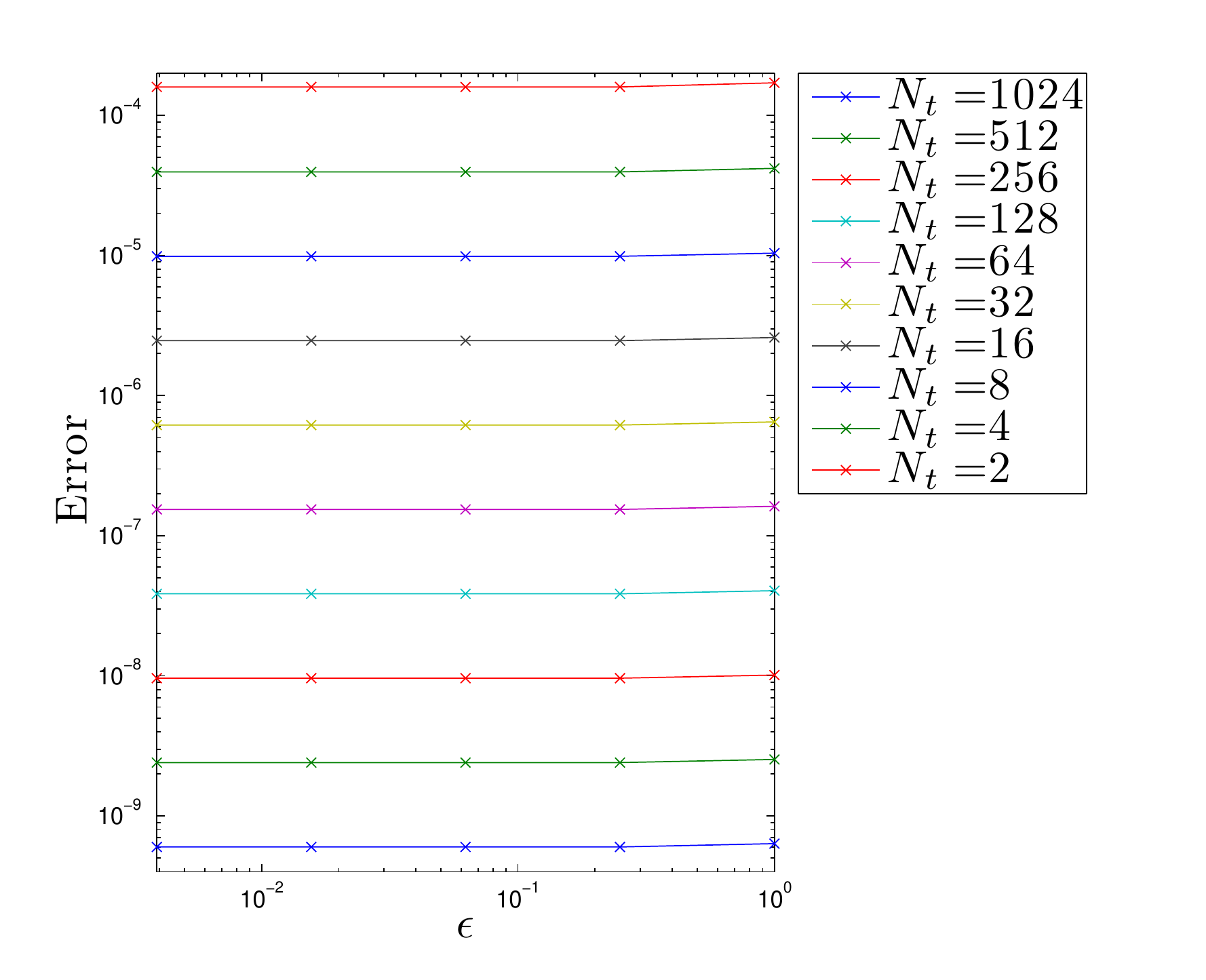}
		\caption{\label{fig:Tf1cvgtpsord2v2NLS} $err_{\rho^\eps} (T_f = 0.1)$ w.r.t $\eps$, $N_x = 2^{9}$} 
	\end{subfigure}
		\hspace{-4cm}
	\caption{\label{fig:Tf1cvgtpsord2-1NLS}Error on the density $\rho^\eps$ for the Strang splitting scheme for \eqref{eq:GPE}  before the caustics: dependence on $\eps$ and on the time step $h$.}
		\end{figure}

\begin{figure}[p]
	\centering
	\hspace{-4cm}
	\begin{subfigure}[t]{0.65\textwidth}
		\centering
		\includegraphics[height=5.cm,width=\textwidth]{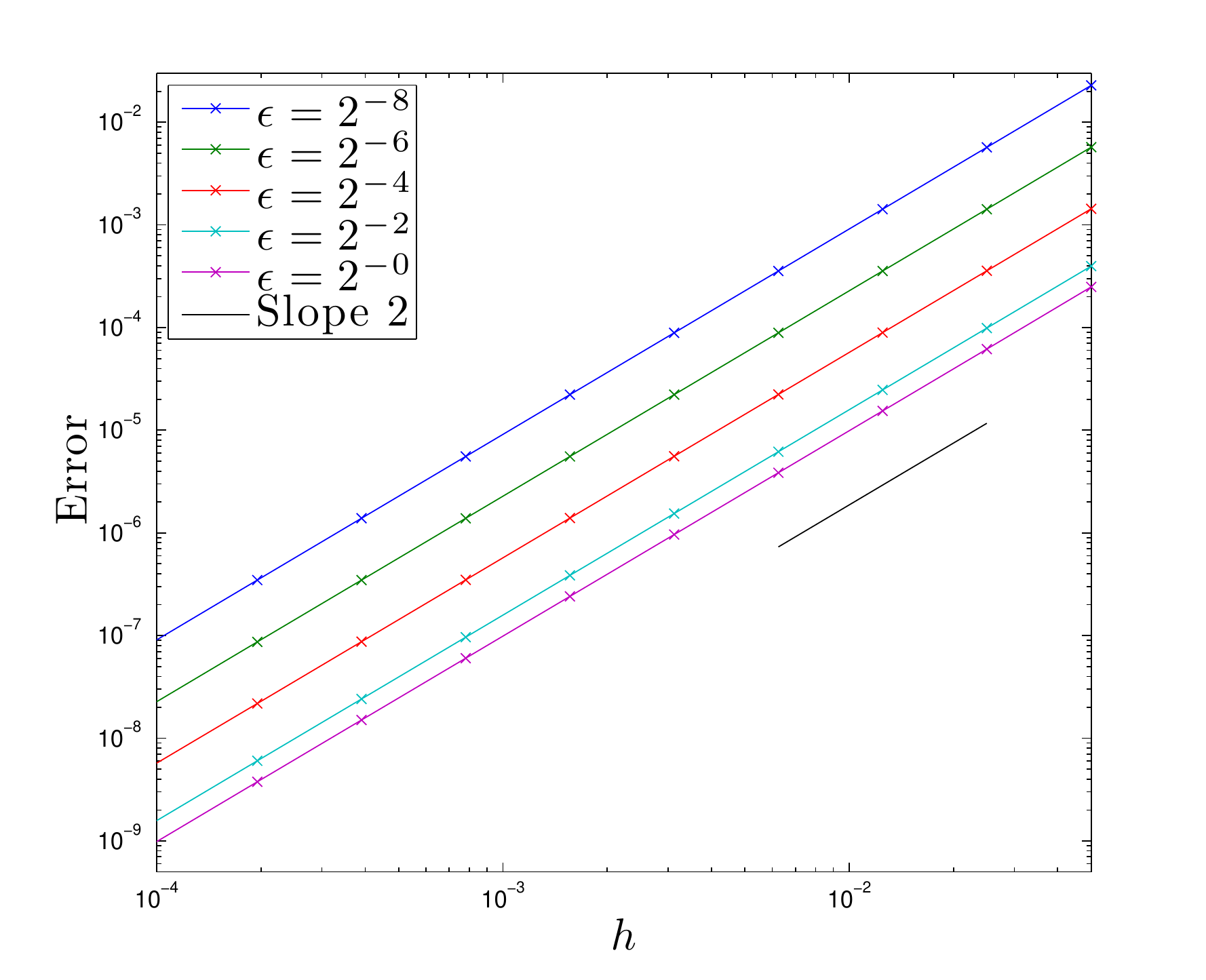}
		\caption{\label{fig:Tf1cvgtpsord2v3NLS} $err_{\pe} (T_f = 0.1)$ w.r.t $h$, $N_x = 2^{9}$}
	\end{subfigure}
	\hspace{-1cm}
	\begin{subfigure}[t]{0.65\textwidth}
		\centering
		\includegraphics[height=5.cm,width=\textwidth]{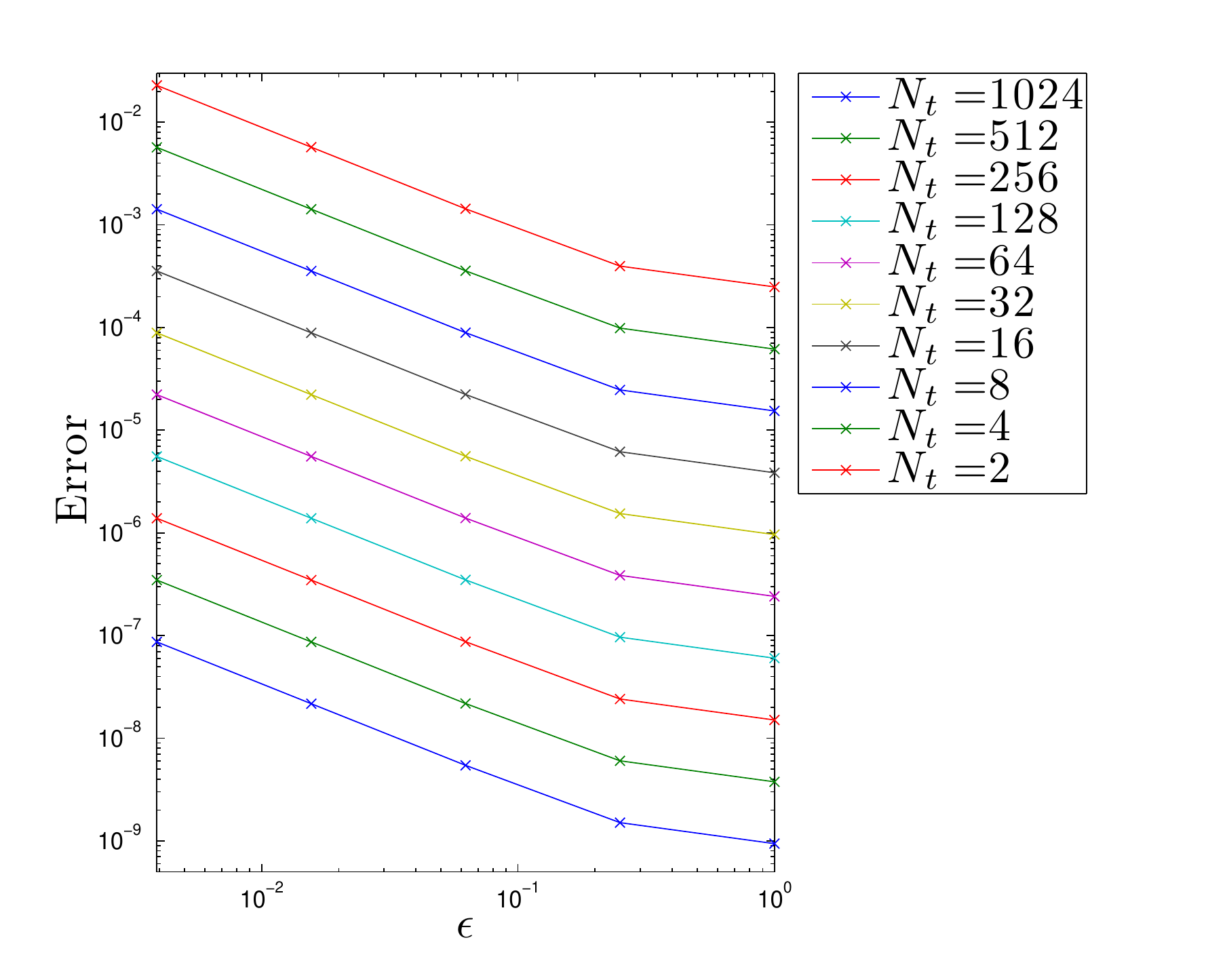}
		\caption{\label{fig:Tf1cvgtpsord2v4NLS} $err_{\pe} (T_f = 0.1)$ w.r.t $\eps$, $N_x = 2^{9}$} 
	\end{subfigure}
	\hspace{-4cm}
	\caption{\label{fig:Tf1cvgtpsord2-2NLS}Error on the wave function $\pe$ for the Strang splitting scheme for \eqref{eq:GPE}  before the caustics: dependence on $\eps$ and on the time step $h$.}
		\end{figure}

\begin{figure}[p]
	\centering
	\hspace{-4cm}
	\begin{subfigure}[t]{0.65\textwidth}
		\centering
		\includegraphics[height=5.cm,width=\textwidth]{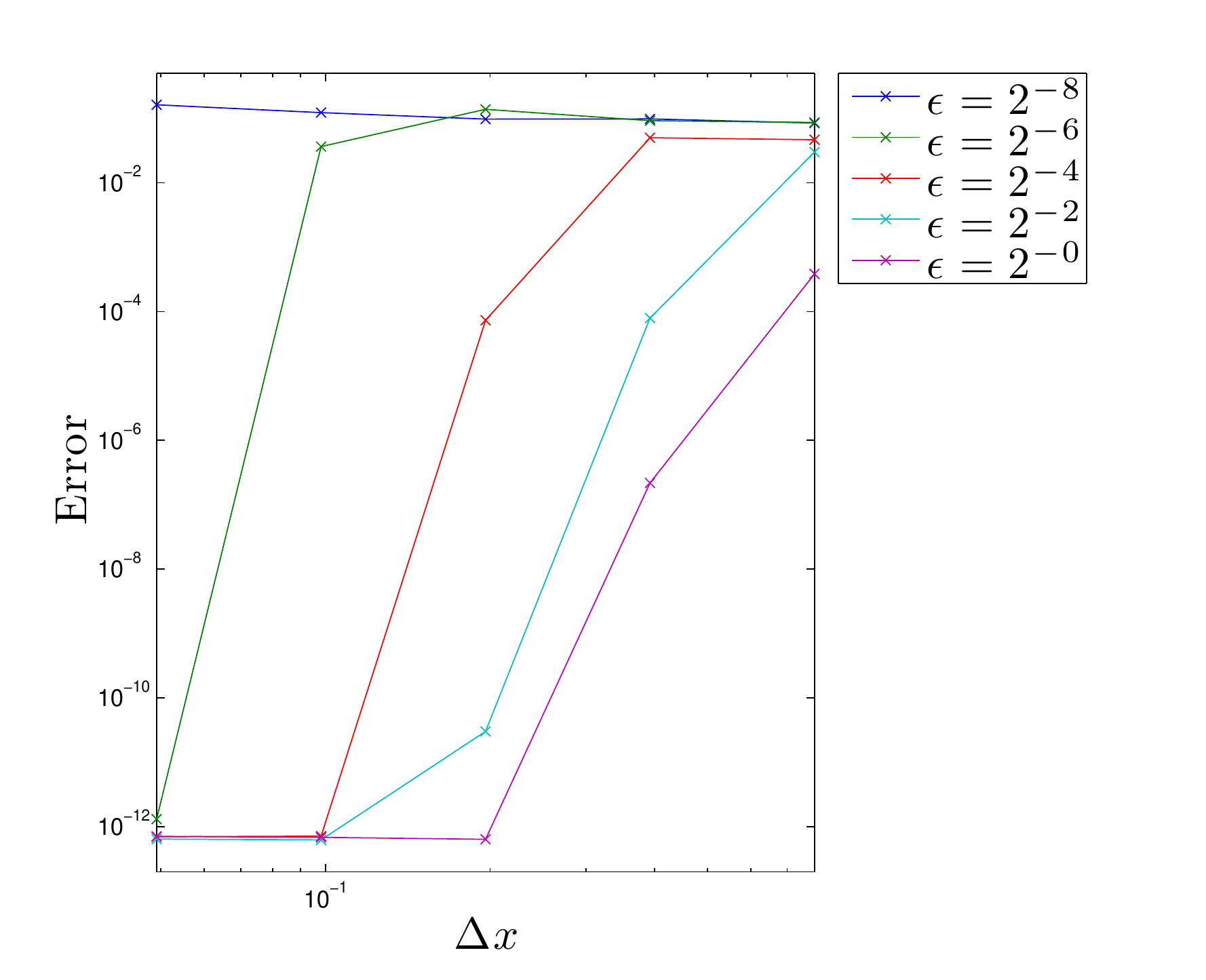}
		\caption{\label{fig:Tf1cvgdxord2v1NLS} $err_{\rho^\eps} (T_f = 0.1)$ w.r.t $\Delta x$, $N_t = 2^{15}$}
	\end{subfigure}
	\hspace{-1cm}
	\begin{subfigure}[t]{0.65\textwidth}
		\centering
		\includegraphics[height=5.cm,width=\textwidth]{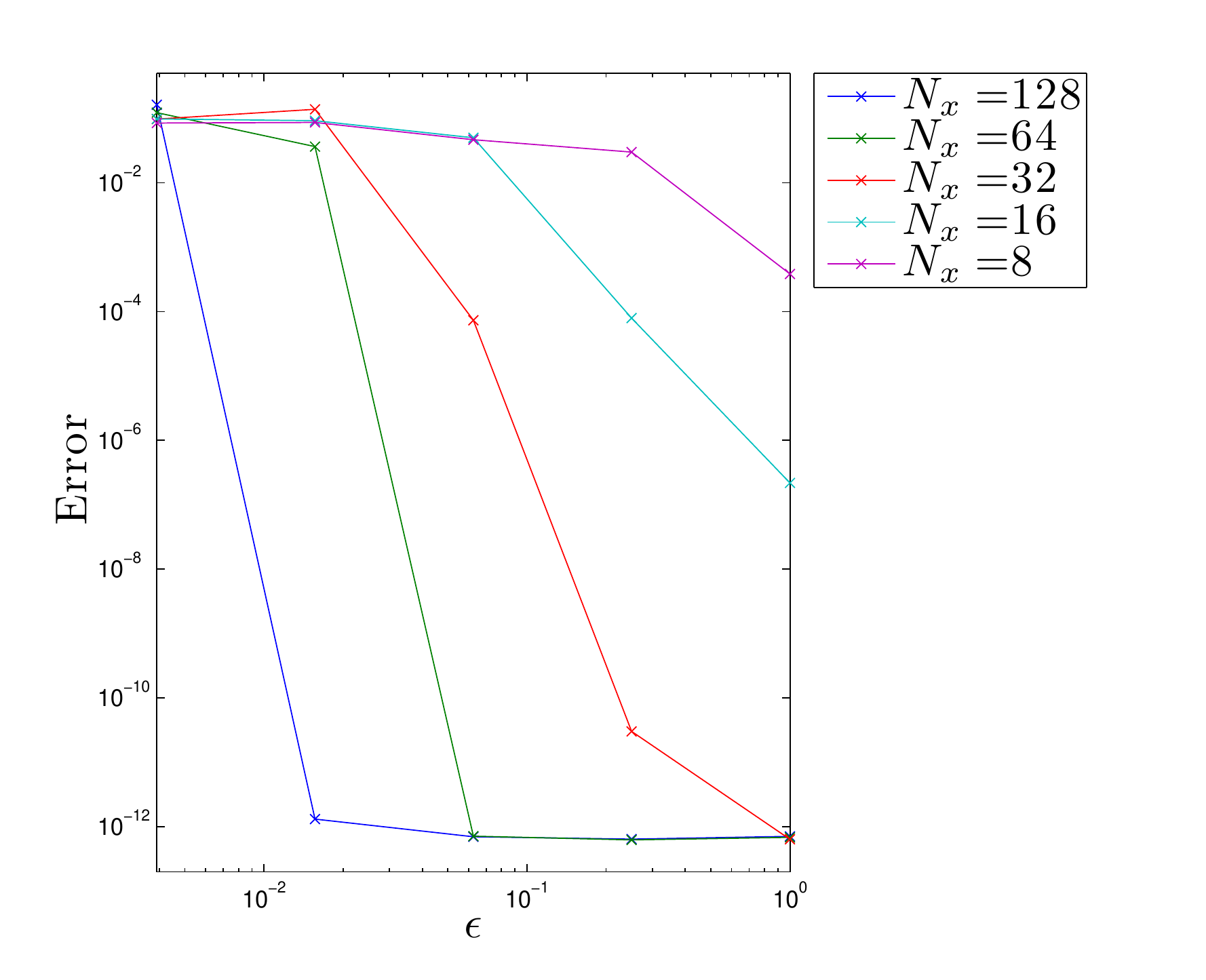}
		\caption{\label{fig:Tf1cvgdxord2v2NLS} $err_{\rho^\eps} (T_f = 0.1)$ w.r.t $\eps$, $N_t = 2^{15}$} 
	\end{subfigure}
	\hspace{-4cm}
	\caption{\label{fig:Tf1cvgdxord2-1NLS}Error on the density $\rho^\eps$ for the Strang splitting scheme for \eqref{eq:GPE}  before the caustics: dependence on $\eps$ and on $\Delta x$.}
		\end{figure}

\begin{figure}[p]
	\centering
	\hspace{-4cm}
	\begin{subfigure}[t]{0.65\textwidth}
		\centering
		\includegraphics[height=5.cm,width=\textwidth]{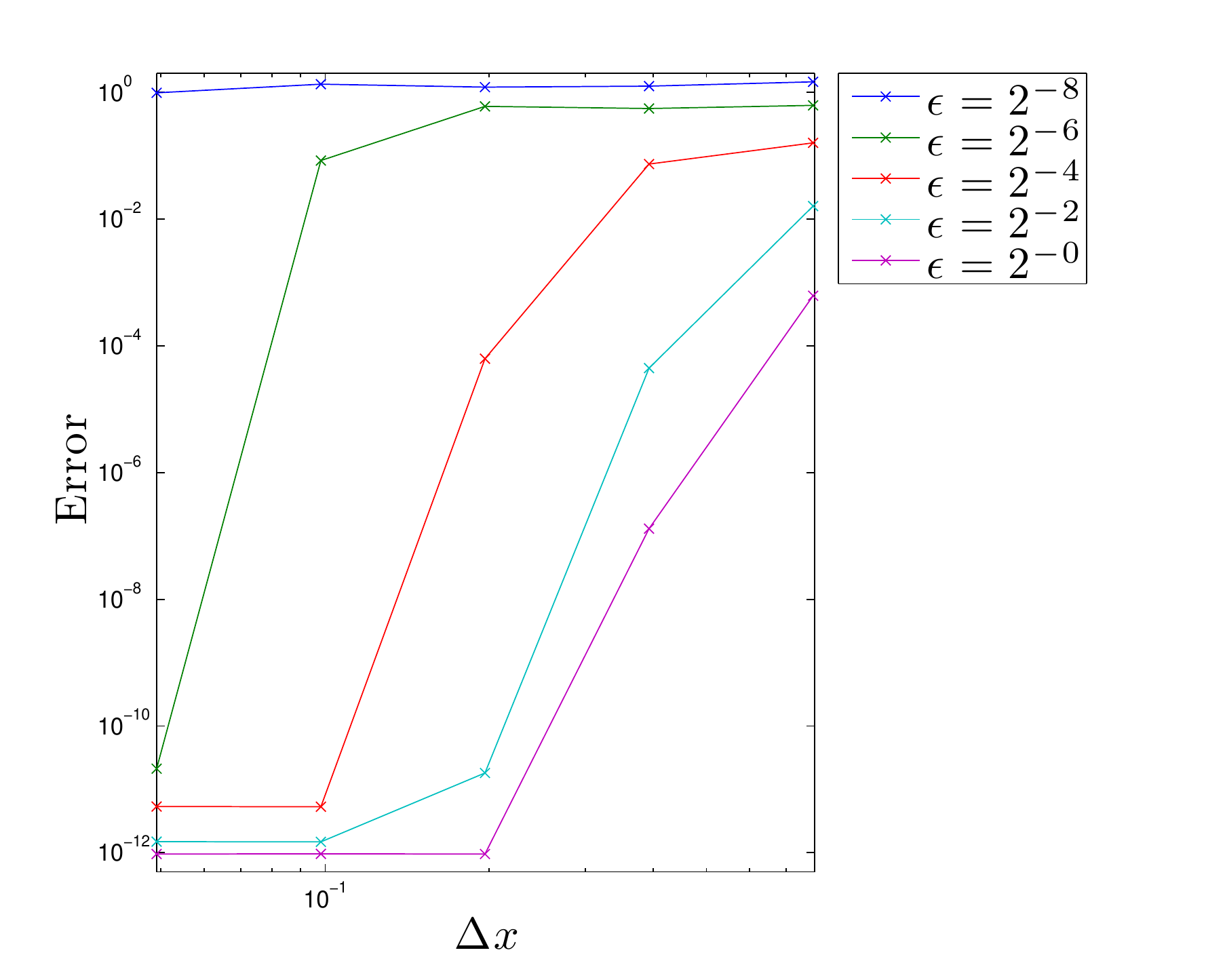}
		\caption{\label{fig:Tf1cvgdxord2v3NLS} $err_{\pe} (T_f = 0.1)$ w.r.t $\Delta x$, $N_t = 2^{15}$}
	\end{subfigure}
	\hspace{-1cm}
	\begin{subfigure}[t]{0.65\textwidth}
		\centering
		\includegraphics[height=5.cm,width=\textwidth]{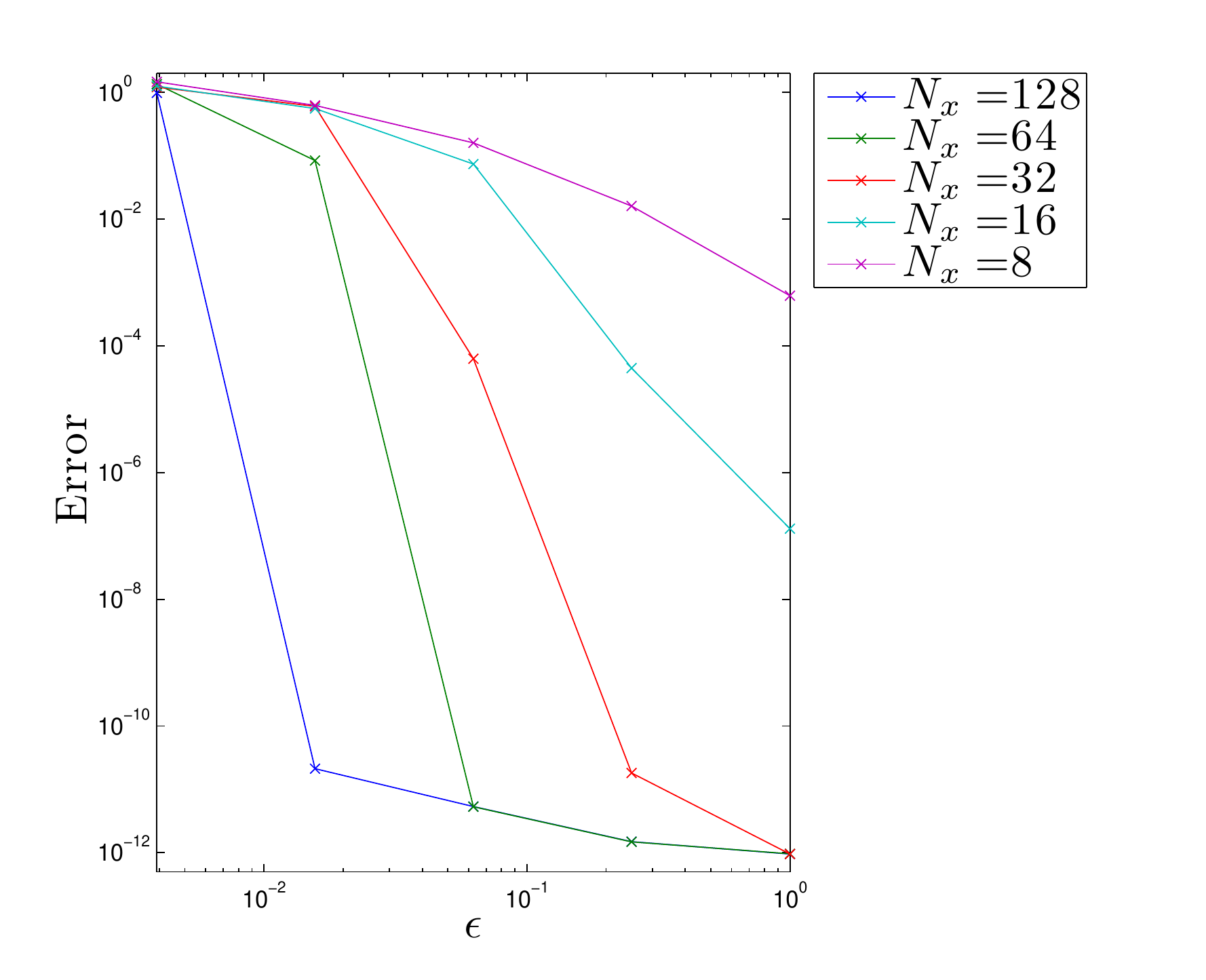}
		\caption{\label{fig:Tf1cvgdxord2v4NLS} $err_{\pe} (T_f = 0.1)$ w.r.t $\eps$, $N_t = 2^{15}$} 
	\end{subfigure}
	\hspace{-4cm}
	\caption{\label{fig:Tf1cvgdxord2-2NLS}Error on the wave function $\pe$ for the Strang splitting scheme for \eqref{eq:GPE}  before the caustics: dependence on $\eps$ and on $\Delta x$.}
	\end{figure}

	Let us now illustrate the behavior of the Strang splitting scheme for \eqref{eq:GPE} at time $T_f = 0.1$ {\em i.e.} before the caustics. 
	On Figures \ref{fig:Tf1cvgtpsord2-1NLS} and \ref{fig:Tf1cvgtpsord2-2NLS}, errors on $\rho^\eps$ and $\pe$ with respect to the time step $h$, for fixed $N_x = 2^9$, are represented and on Figures \ref{fig:Tf1cvgdxord2-1NLS} to \ref{fig:Tf1cvgdxord2-2NLS}, errors with respect to $\Delta x$, for fixed $N_t = h/T_f = 2^{15}$, are represented. Regarding the observable $\rho^\eps = |\Psi^{\eps}|^2$, \ref{fig:Tf1cvgtpsord2v1NLS}, \ref{fig:Tf1cvgtpsord2v2NLS}, \ref{fig:Tf1cvgdxord2v1NLS} and \ref{fig:Tf1cvgdxord2v2NLS} corroborate the fact that the error behaves as
	\[
		\calO \(h^2 + C_{\eps,N}\Delta x^N\)
	\]
	where $N>0$ and $C_{\eps,N}\to +\infty$ as $\eps\to 0$ \cite{MR1880116,MR2047194,MR3138106}. This is in agreement with the results obtained by Carles \cite{MR3138106} in the weakly nonlinear case before the caustics; however, our simulations suggest that this behavior persists in the supercritical case. If we observe the wave function, the situation is completely different: the Strang splitting scheme is not UA any more when $h\to 0$. Figures \ref{fig:Tf1cvgtpsord2v3NLS}, \ref{fig:Tf1cvgtpsord2v4NLS}, \ref{fig:Tf1cvgdxord2v3NLS} and \ref{fig:Tf1cvgdxord2v4NLS} indeed suggest that the error of $\Psi^\eps$ behaves like 
	\[
		\calO \(\frac{h^2}{\eps} + C_{\eps,N}\Delta x^N\)
	\]
	where $N>0$ and $C_{\eps,N}\to +\infty$ as $\eps\to 0$.
	
\begin{figure}[p]
	\centering
	\hspace{-4cm}
	\begin{subfigure}[t]{0.6\textwidth}
		\centering
		\includegraphics[height=5.cm,width=\textwidth]{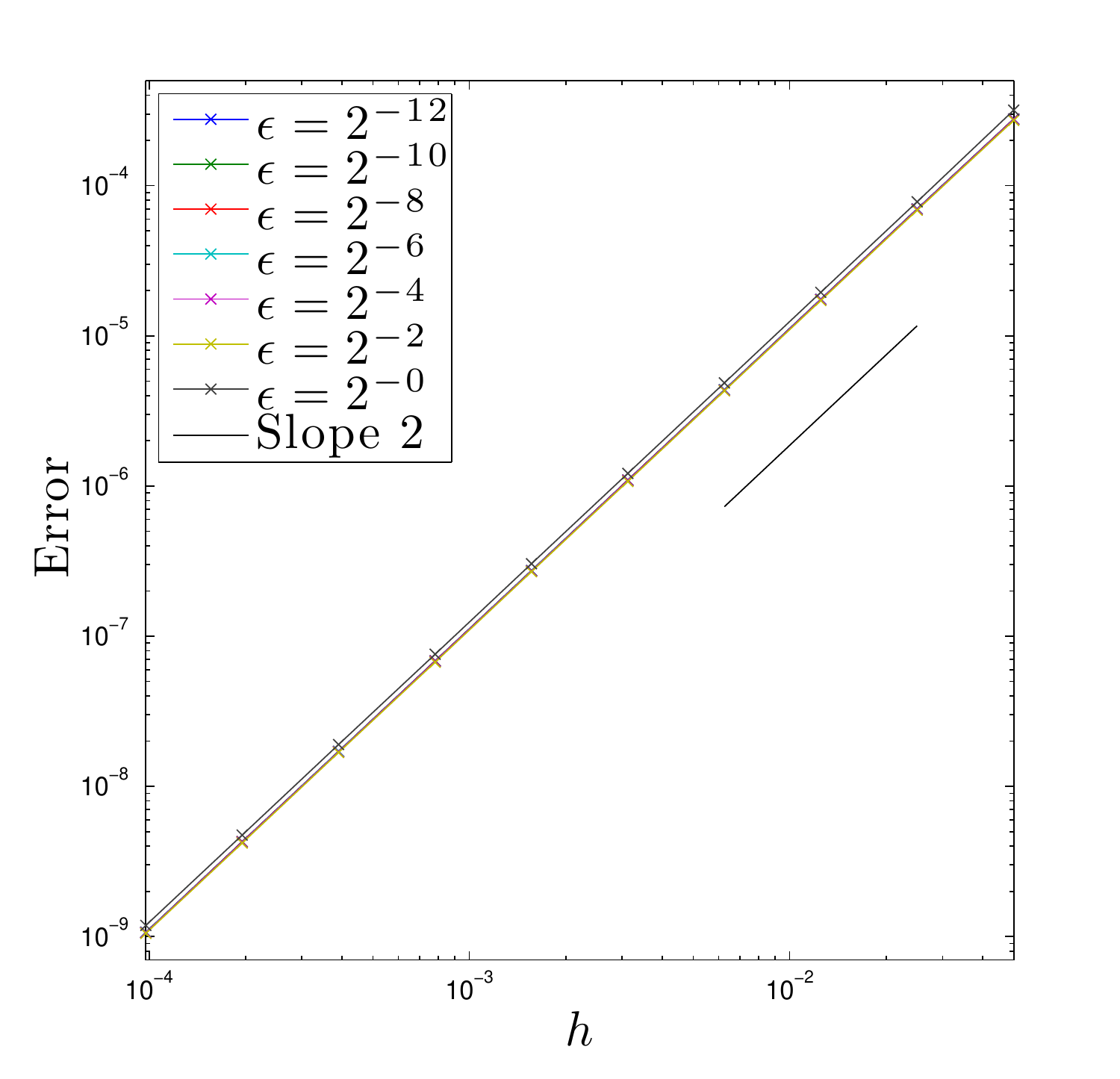}
		\caption{\label{fig:Tf1cvgtpsord2v1} $err_{\rho^\eps} (T_f = 0.1)$ w.r.t $h$, $N_x = 2^{7}$}
	\end{subfigure}
	\hspace{-1cm}
	\begin{subfigure}[t]{0.7\textwidth}
		\centering
		\includegraphics[height=5.cm,width=\textwidth]{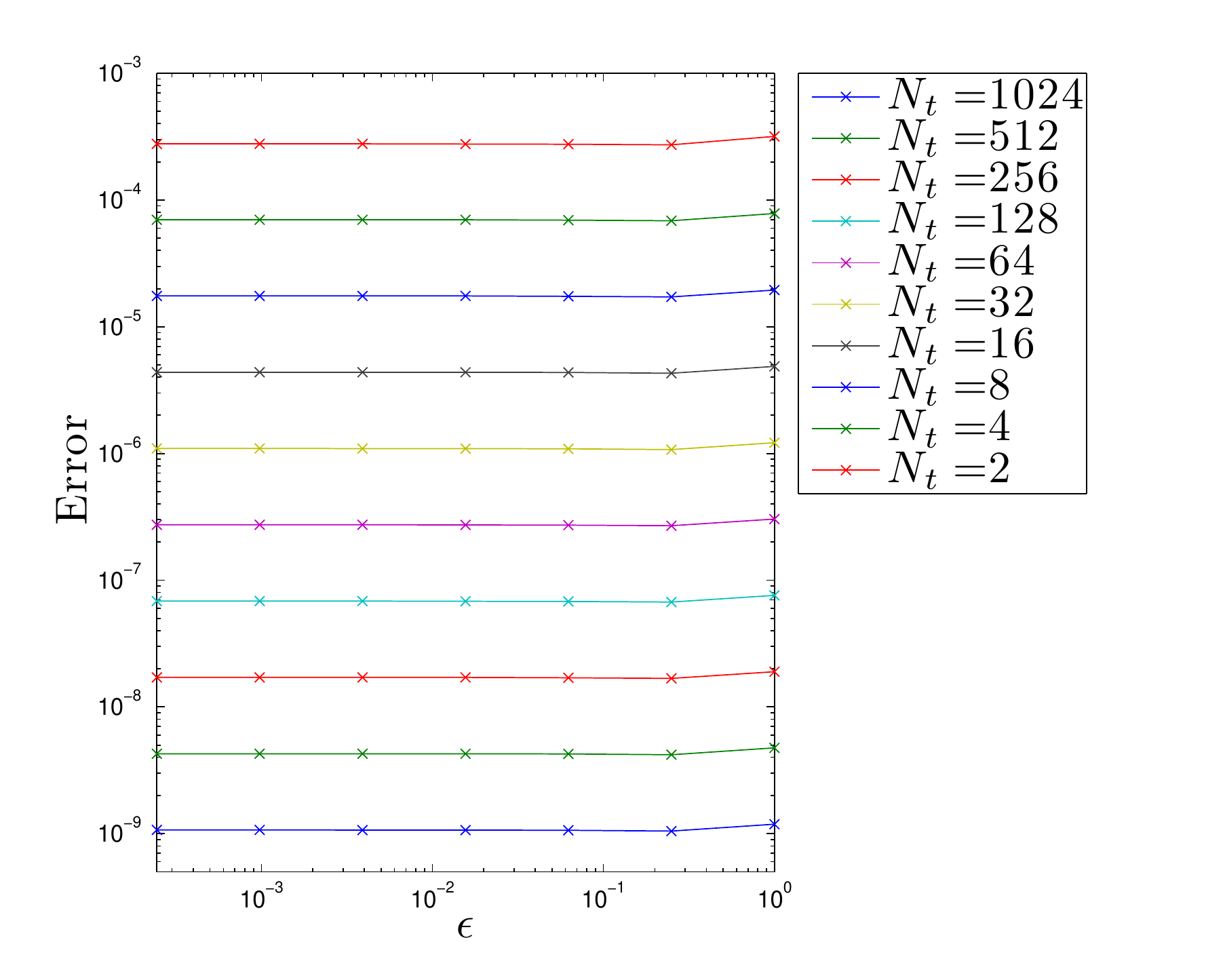}
		\caption{\label{fig:Tf1cvgtpsord2v2} $err_{\rho^\eps} (T_f = 0.1)$ w.r.t $\eps$, $N_x = 2^{7}$} 
	\end{subfigure}
	\hspace{-4cm}
	\caption{\label{fig:Tf1cvgtpsord2-1} Error on the density $\rho^\eps$ for the splitting scheme \eqref{scheme2} of order $2$ before the caustics: dependence on $\eps$ and on $h$.}
		\end{figure}

\begin{figure}[p]
	\centering
	\hspace{-4cm}
	\begin{subfigure}[t]{0.6\textwidth}
		\centering
		\includegraphics[height=5.cm,width=\textwidth]{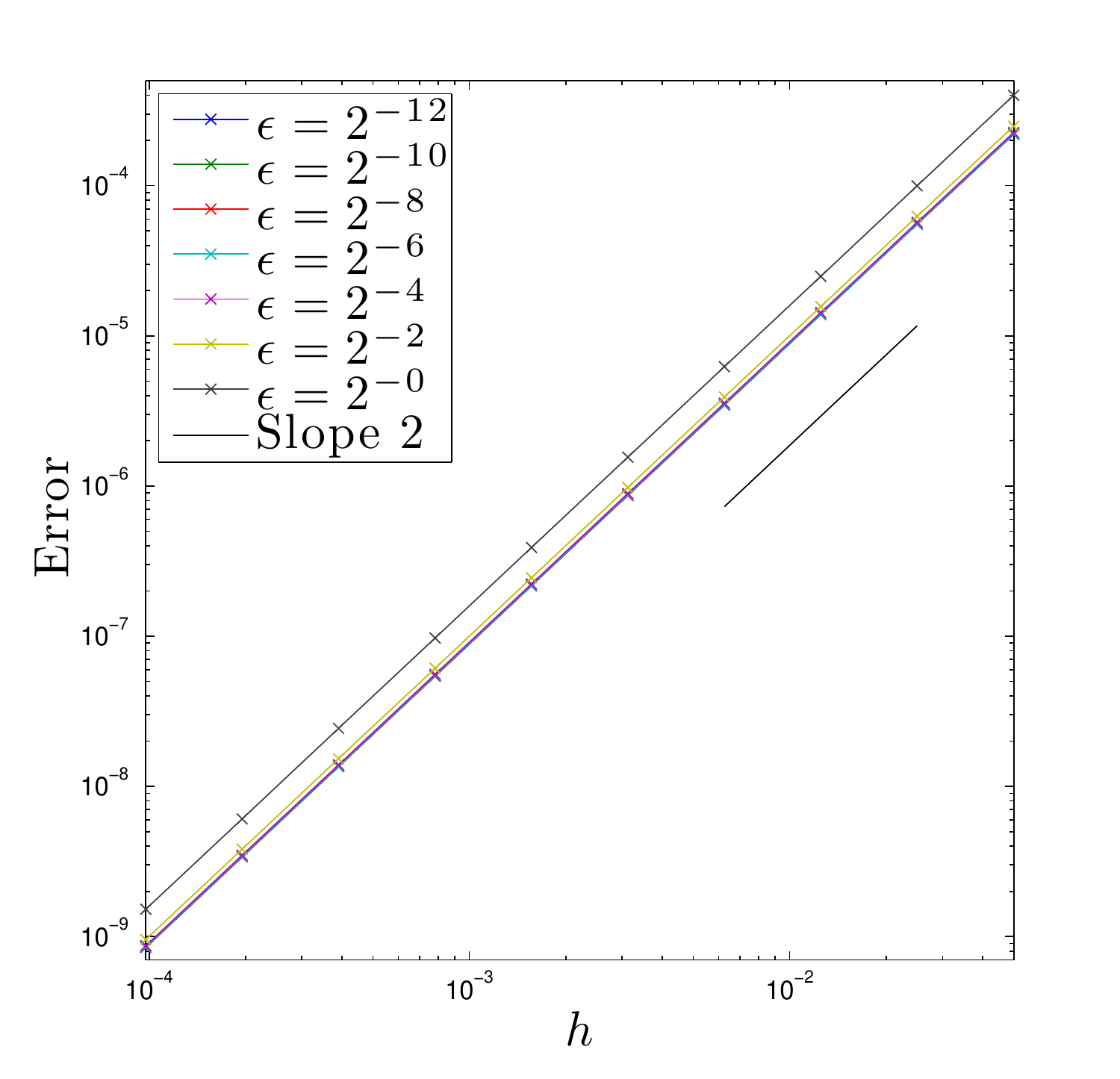}
		\caption{\label{fig:Tf1cvgtpsord2v3} $err_{(\se,\ae)} (T_f = 0.1)$ w.r.t $h$, $N_x = 2^{7}$}
	\end{subfigure}
	\hspace{-1cm}
	\begin{subfigure}[t]{0.7\textwidth}
		\centering
		\includegraphics[height=5.cm,width=\textwidth]{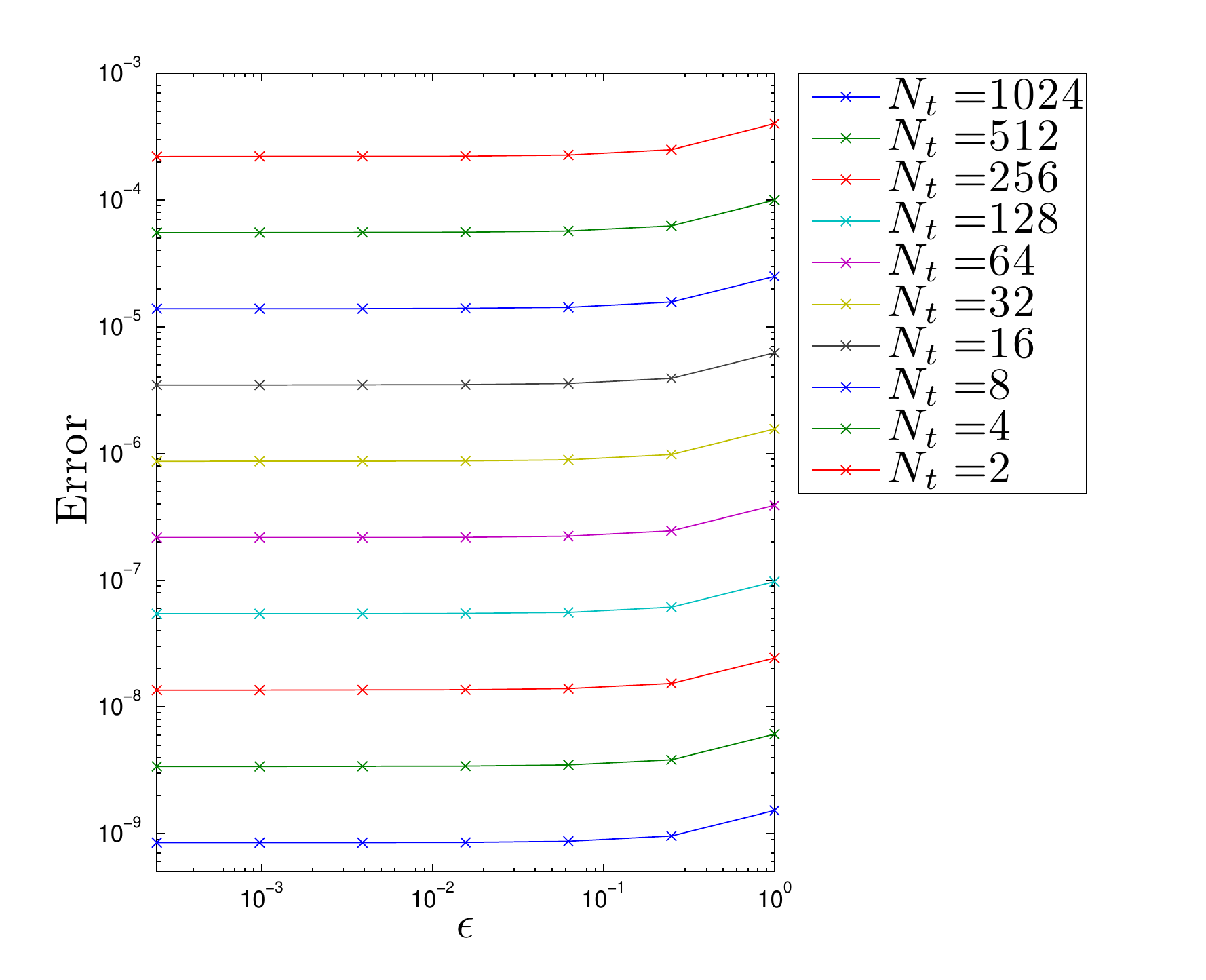}
		\caption{\label{fig:Tf1cvgtpsord2v4} $err_{(\se,\ae)} (T_f = 0.1)$ w.r.t $\eps$, $N_x = 2^{7}$} 
	\end{subfigure}
	\hspace{-4cm}
	\caption{\label{fig:Tf1cvgtpsord2-2} Error on $(\se,\ae)$ for the splitting scheme \eqref{scheme2} of order $2$ before the caustics: dependence on $\eps$ and on $h$.}
		\end{figure}

\begin{figure}[p]
	\centering
	\hspace{-4cm}
	\begin{subfigure}[t]{0.65\textwidth}
		\centering
		\includegraphics[height=5.cm,width=\textwidth]{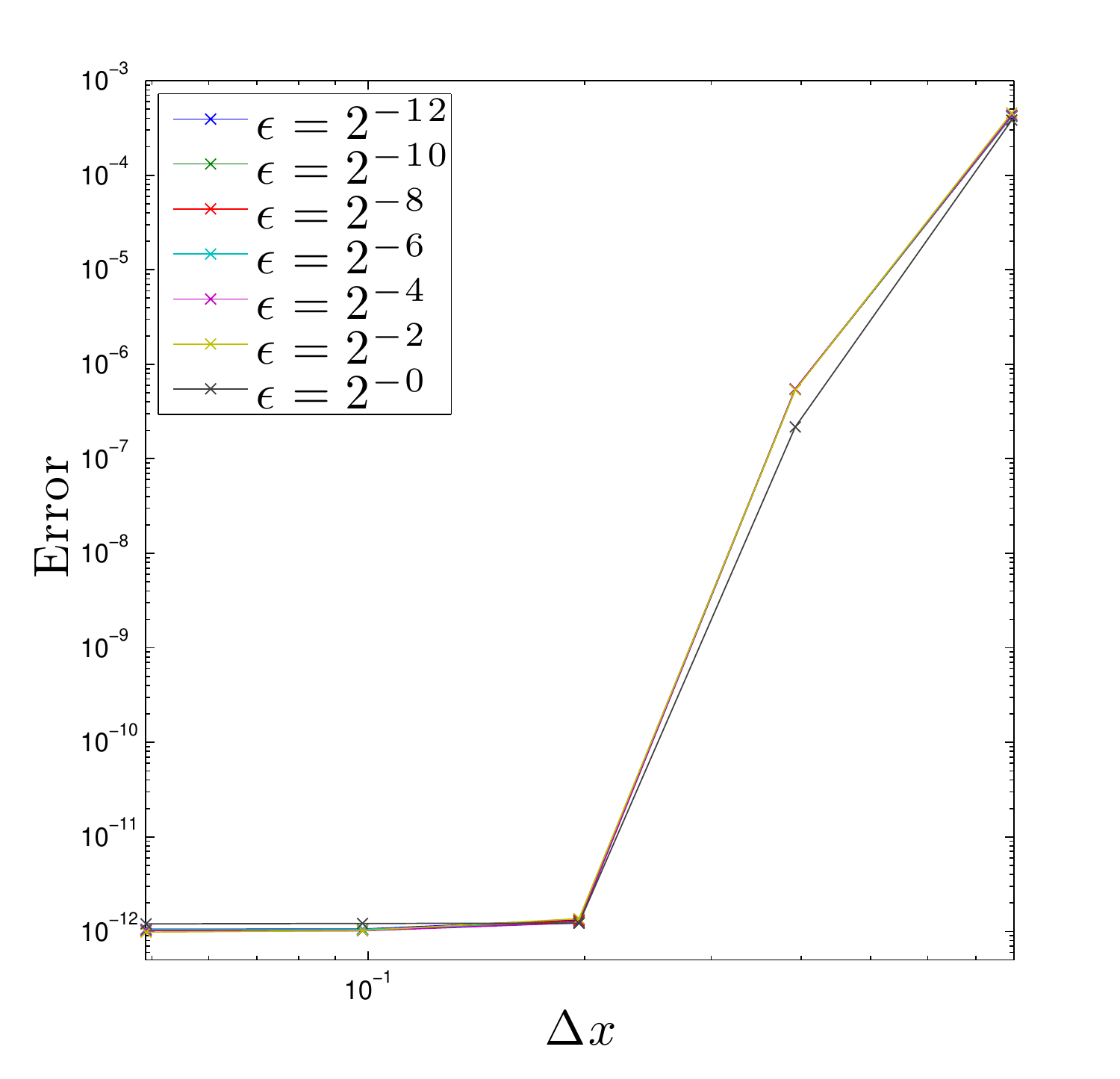}
		\caption{\label{fig:Tf1cvgdxord2v1} $err_{\rho^\eps} (T_f = 0.1)$ w.r.t $\Delta x$, $N_t = 2^{15}$}
	\end{subfigure}
	\hspace{-1cm}
	\begin{subfigure}[t]{0.65\textwidth}
		\centering
		\includegraphics[height=5.cm,width=\textwidth]{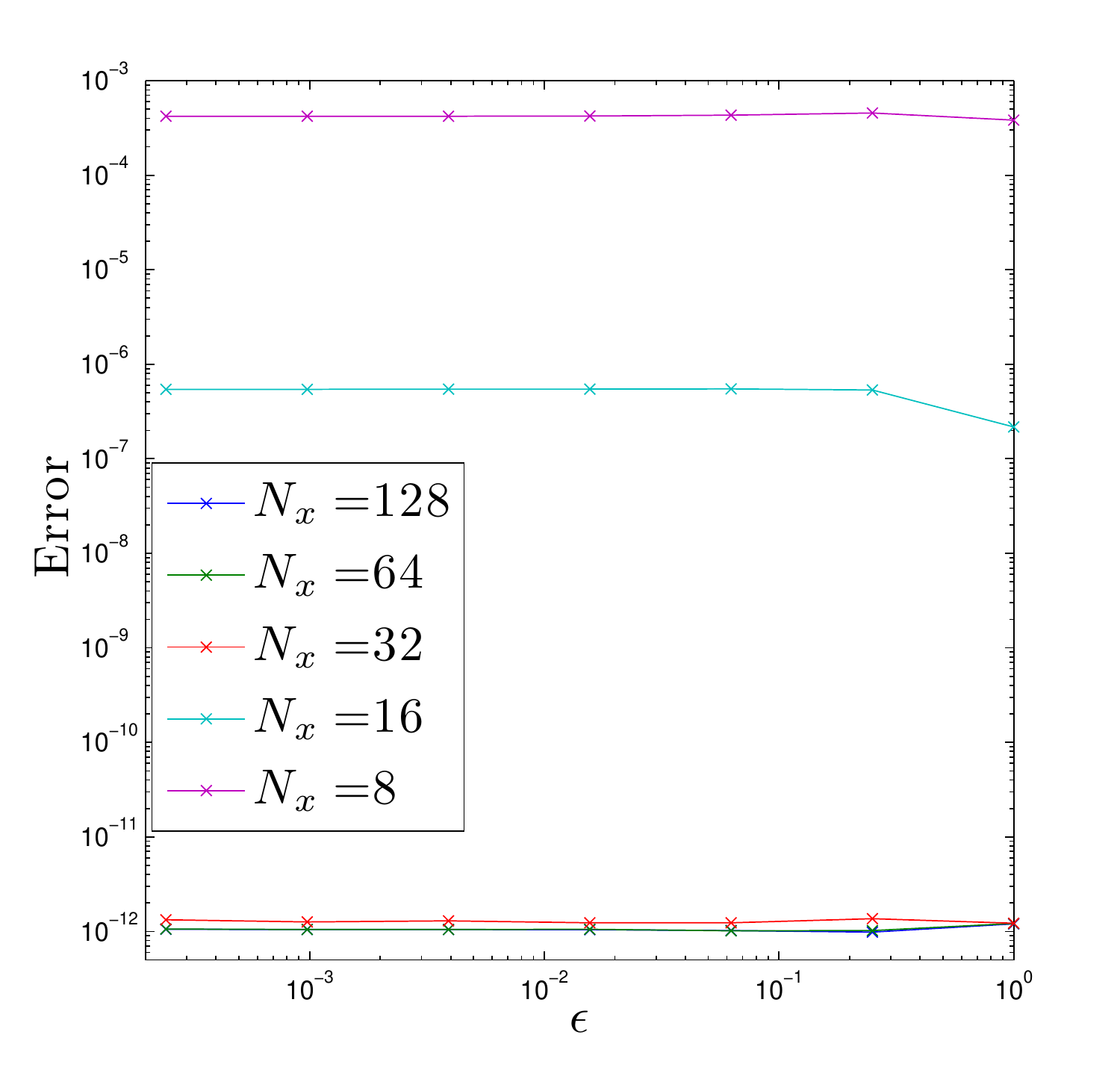}
		\caption{\label{fig:Tf1cvgdxord2v2} $err_{\rho^\eps} (T_f = 0.1)$ w.r.t $\eps$, $N_t = 2^{15}$} 
	\end{subfigure}
	\hspace{-4cm}
		\caption{\label{fig:Tf1cvgdxord2-1} Error on the density $\rho^\eps$ for the splitting scheme \eqref{scheme2} of order $2$ before the caustics: dependence on $\eps$ and on $\Delta x$.}

		\end{figure}

\begin{figure}[p]
	\centering
	\hspace{-4cm}
	\begin{subfigure}[t]{0.65\textwidth}
		\centering
		\includegraphics[height=5.cm,width=\textwidth]{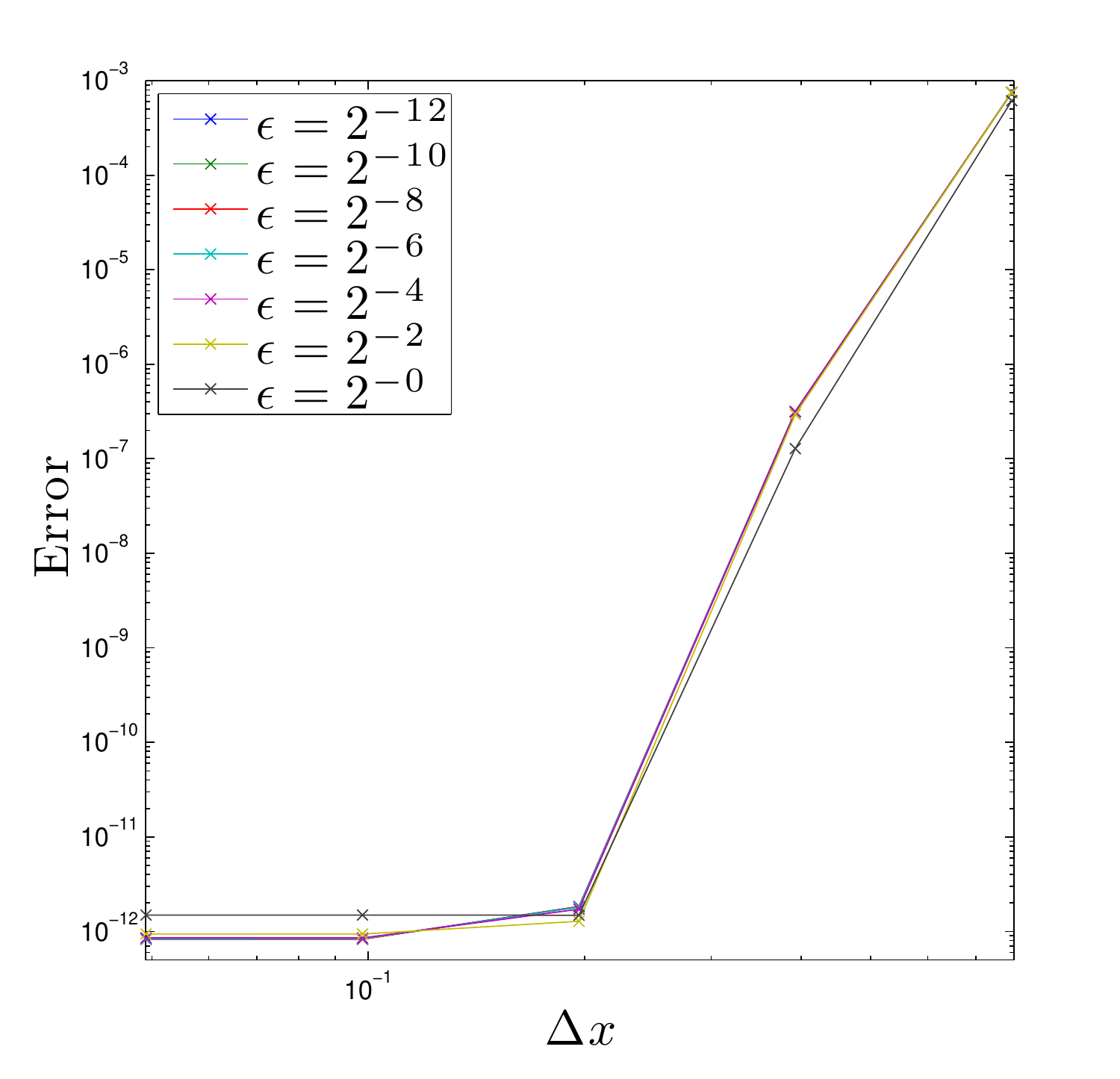}
		\caption{\label{fig:Tf1cvgdxord2v3} $err_{(\se,\ae)} (T_f = 0.1)$ w.r.t $\Delta x$, $N_t = 2^{15}$}
	\end{subfigure}
	\hspace{-1cm}
	\begin{subfigure}[t]{0.65\textwidth}
		\centering
		\includegraphics[height=5.cm,width=\textwidth]{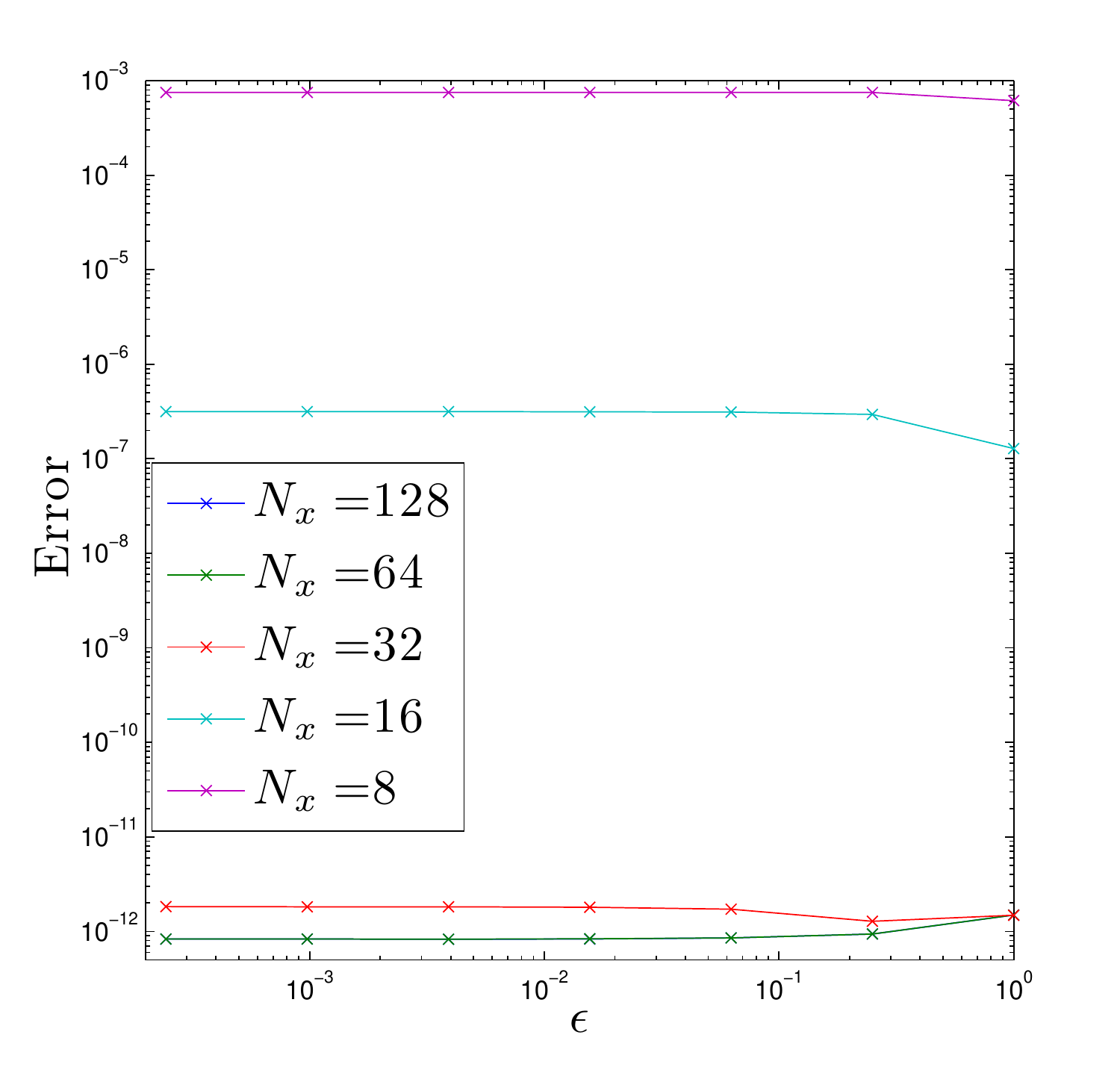}
		\caption{\label{fig:Tf1cvgdxord2v4} $err_{(\se,\ae)} (T_f = 0.1)$ w.r.t $\eps$, $N_t = 2^{15}$} 
	\end{subfigure}
	\hspace{-4cm}
	\caption{\label{fig:Tf1cvgdxord2-2} Error on $(\se,\ae)$ for the splitting scheme \eqref{scheme2} of order $2$ before the caustics: dependence on $\eps$ and on $\Delta x$.}		\end{figure}

	Let us now focus on the experiments performed with our second and fourth-order methods, in the same situation. We start with the second-order scheme \eqref{scheme2}. Figures \ref{fig:Tf1cvgtpsord2-1} and \ref{fig:Tf1cvgtpsord2-2} represent the errors on $\rho^\eps$ and $(\se,\ae)$ {\em w.r.t.}{\@} the time step $h$ for a fixed $N_x = 2^7$. Figures \ref{fig:Tf1cvgdxord2-1} and \ref{fig:Tf1cvgdxord2-2} represent the errors {\em w.r.t.}\@ $\Delta x$ for fixed $N_t = h/T_f = 2^{13}$. All these figures illustrate the fact that our scheme is UA with respect to $\eps$, for the quadratic observables as well as for the whole unknown $(\se,\ae)$ itself. Figures \ref{fig:Tf1cvgtpsord2-1} and \ref{fig:Tf1cvgtpsord2-2} show that \eqref{scheme2} is uniformly of order 2 in time, whereas Figures  \ref{fig:Tf1cvgdxord2-1} and \ref{fig:Tf1cvgdxord2-2} show that the convergence is uniformly spectral in space. 
	
\begin{figure}[p]
	\centering	
		\hspace{-5cm}
		\begin{subfigure}[t]{0.6\textwidth}
		\centering
		\includegraphics[height=5.cm,width=\textwidth]{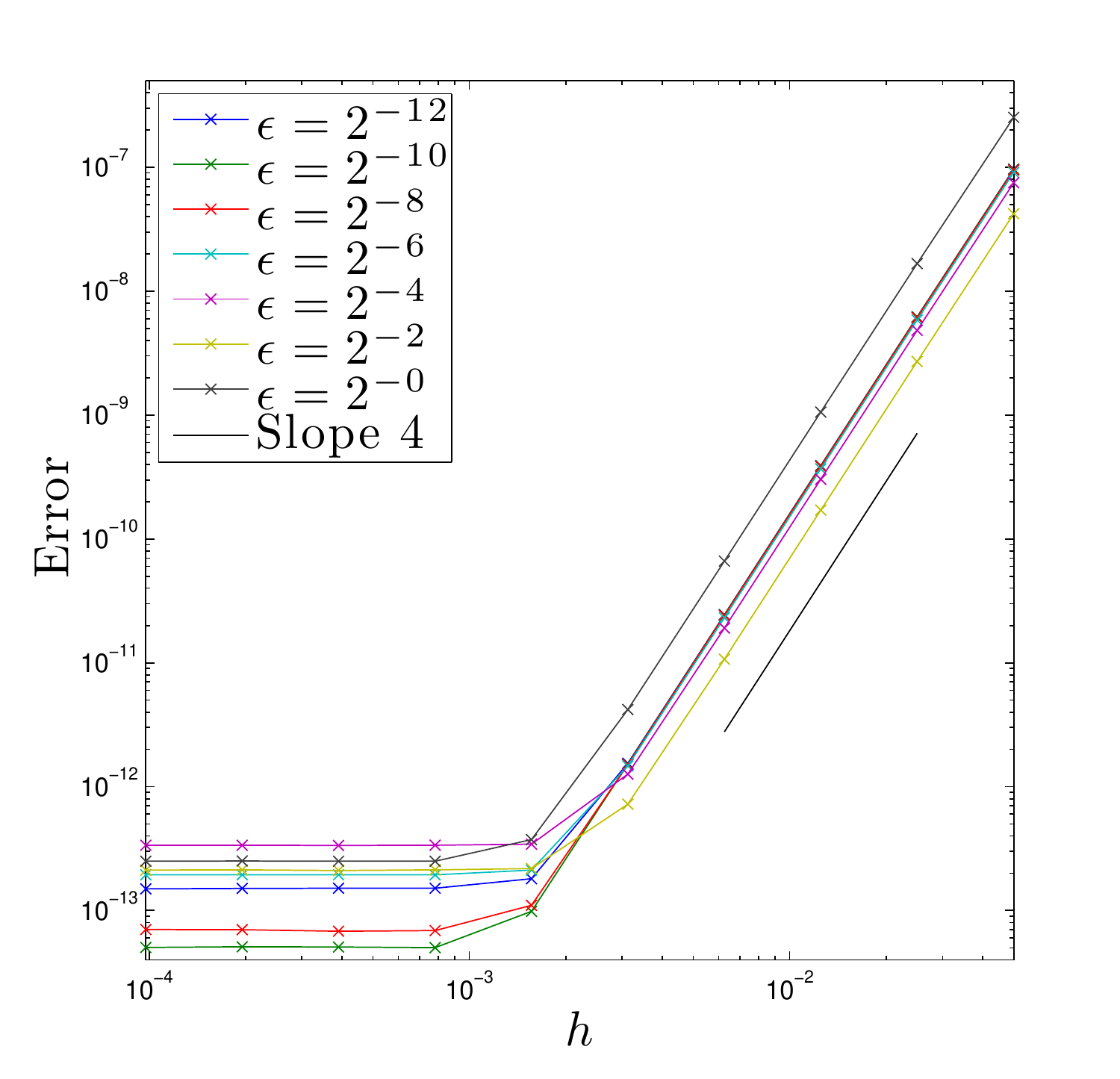}
		\caption{\label{fig:Tf1cvgtpsord4v1} $err_{\rho^\eps} (T_f = 0.1)$ w.r.t $h$, $N_x = 2^7$ } 
		\end{subfigure}
		\hspace{-1cm}
		\begin{subfigure}[t]{0.71\textwidth}
		\centering
		\includegraphics[height=5.cm,width=\textwidth]{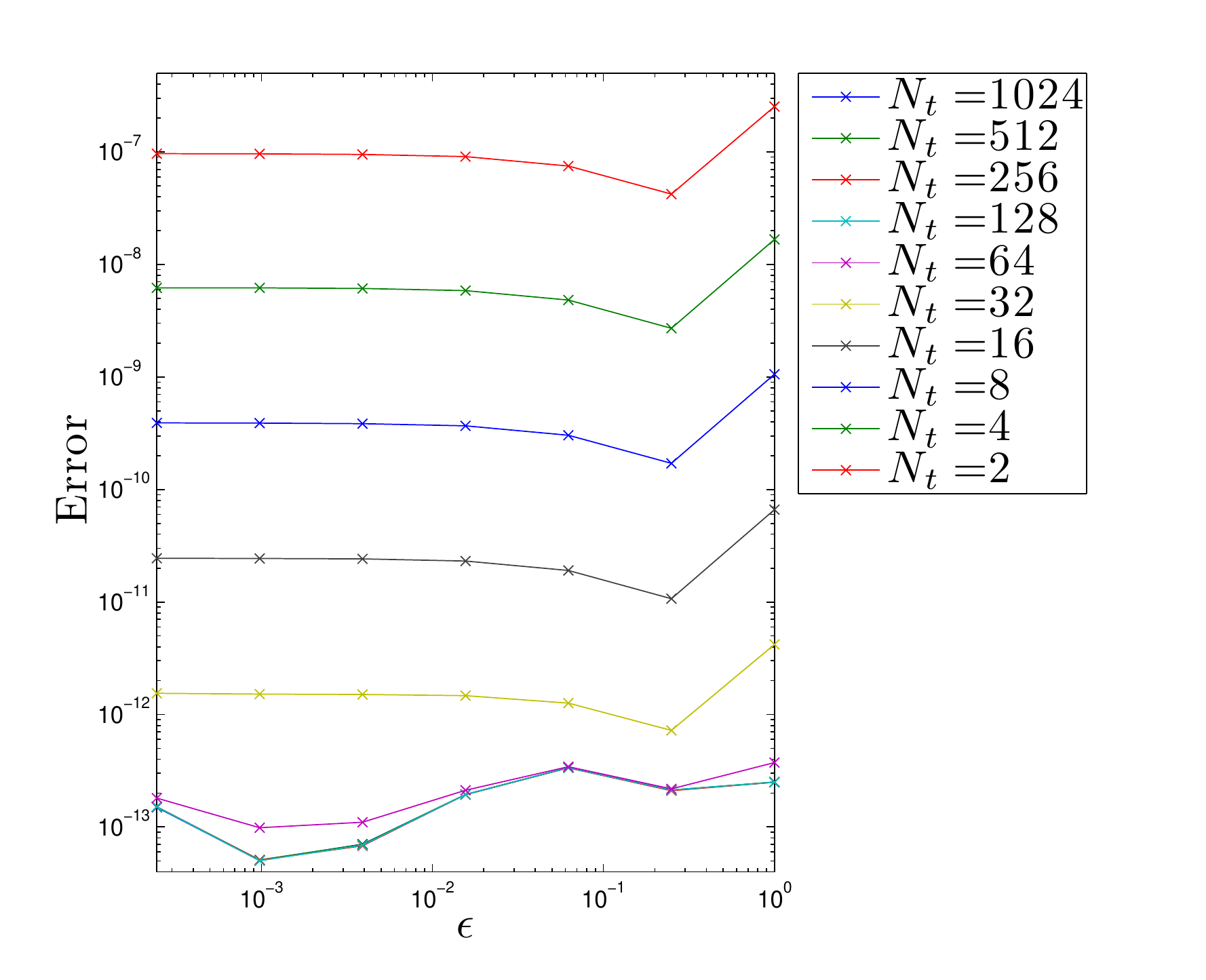}
		\caption{\label{fig:Tf1cvgtpsord4v2} $err_{\rho^\eps} (T_f = 0.1)$ w.r.t $\eps$, $N_x = 2^7$} 
		\end{subfigure}
		\hspace{-6cm}
		\caption{\label{fig:Tf1cvgtpsord4-1}Error on the density $\rho^\eps$ for the splitting scheme \eqref{scheme4} of order $4$ before the caustics: dependence on $\eps$ and on $h$.}

		\end{figure}

\begin{figure}[p]
		\centering	
		\hspace{-5cm}
		\begin{subfigure}[t]{0.6\textwidth}
		\centering
		\includegraphics[height=5.cm,width=\textwidth]{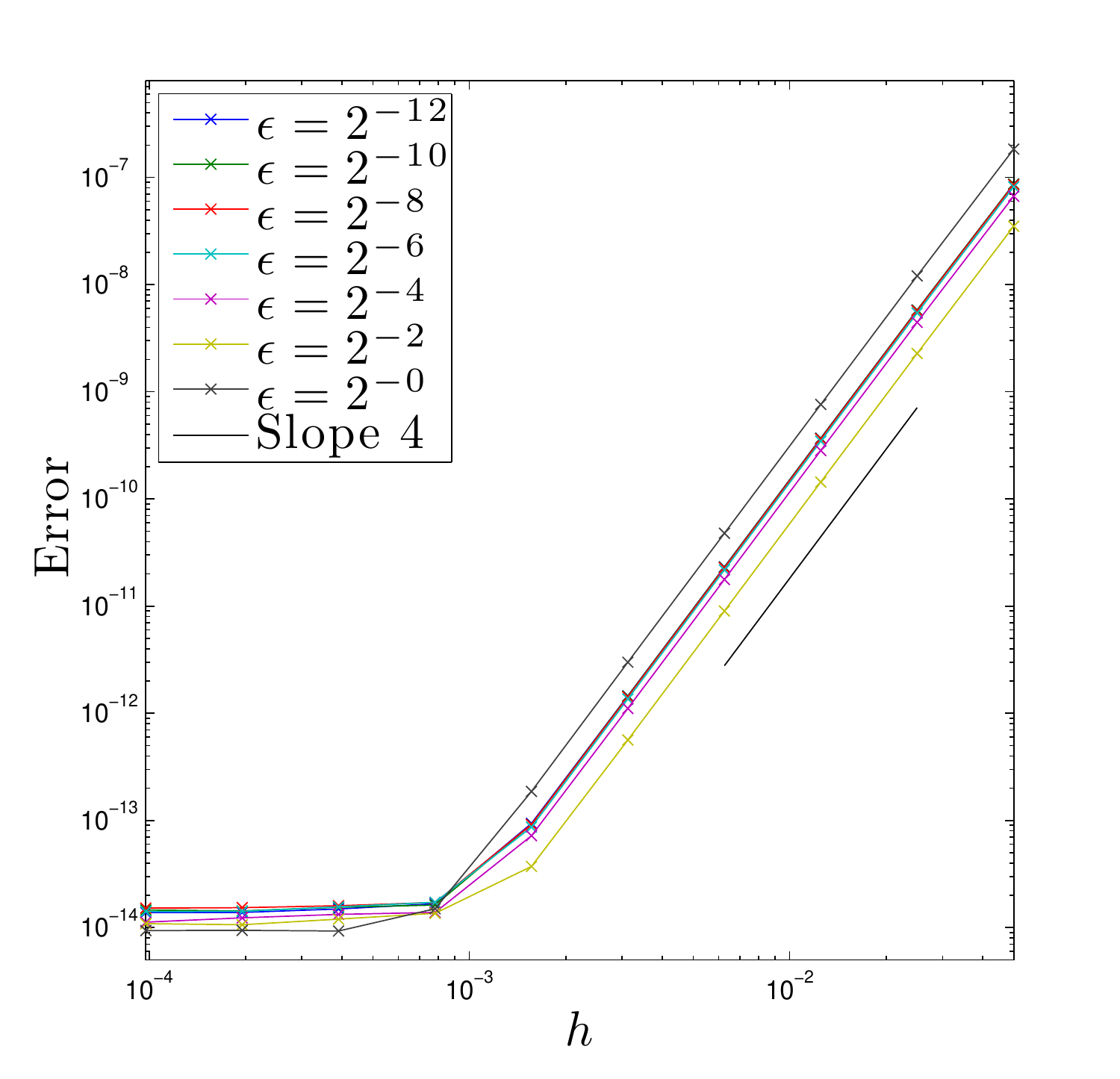}
		\caption{\label{fig:Tf1cvgtpsord4v3} $err_{(\se,\ae)} (T_f = 0.1)$ w.r.t $h$, $N_x = 2^7$ } 
		\end{subfigure}
		\hspace{-1cm}
		\begin{subfigure}[t]{0.71\textwidth}
		\centering
		\includegraphics[height=5.cm,width=\textwidth]{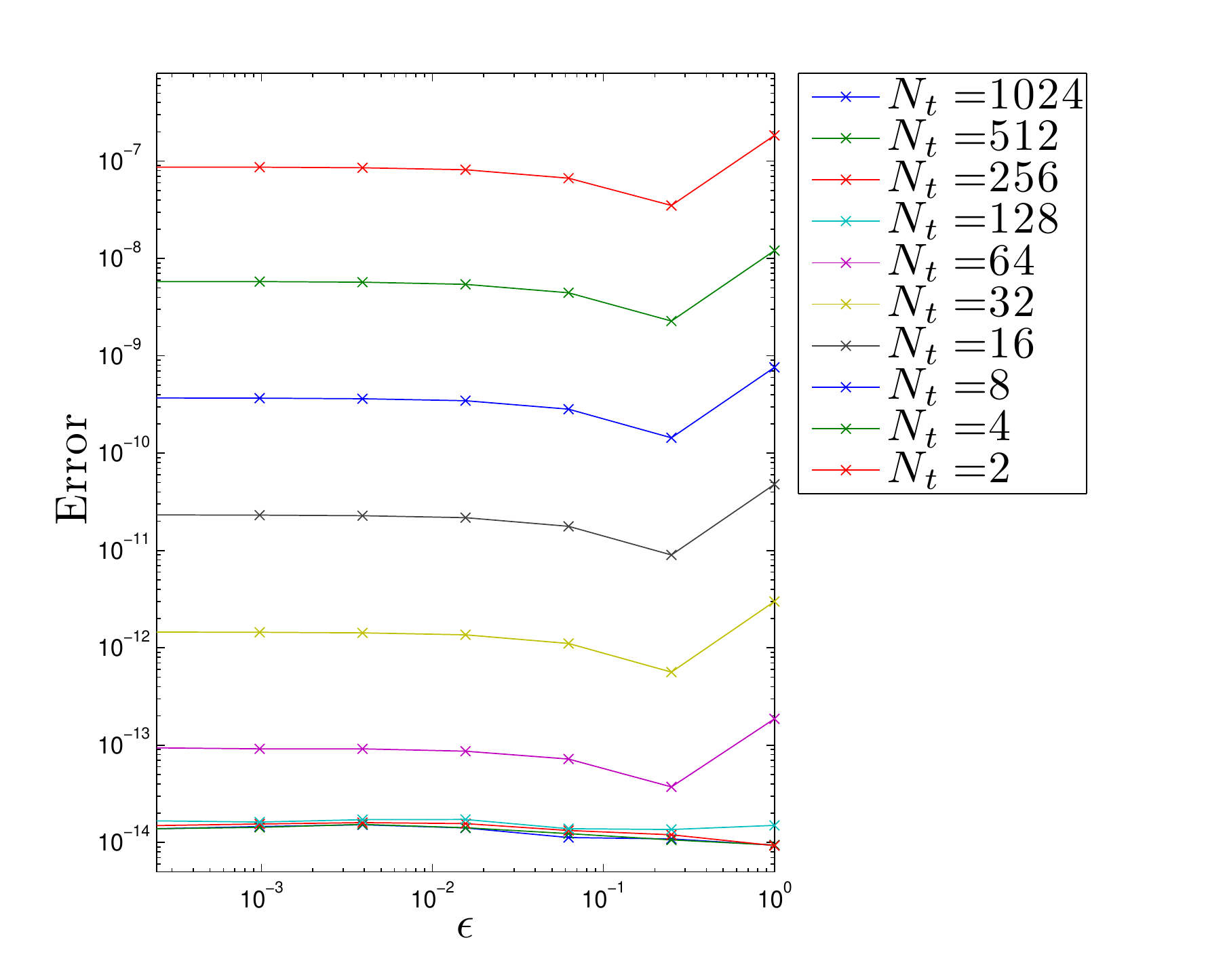}
		\caption{\label{fig:Tf1cvgtpsord4v4} $err_{(\se,\ae)} (T_f = 0.1)$ w.r.t $\eps$, $N_x = 2^7$} 
		\end{subfigure}
		\hspace{-6cm}
		\caption{ \label{fig:Tf1cvgtpsord4-2}Error on $(\se,\ae)$ for the splitting scheme \eqref{scheme4} of order $4$ before the caustics: dependence on $\eps$ and on $h$.}

		\end{figure}

\begin{figure}[p]
		\centering	
		\hspace{-5cm}
		\begin{subfigure}[t]{0.6\textwidth}
		\centering
		\includegraphics[height=5.cm,width=\textwidth]{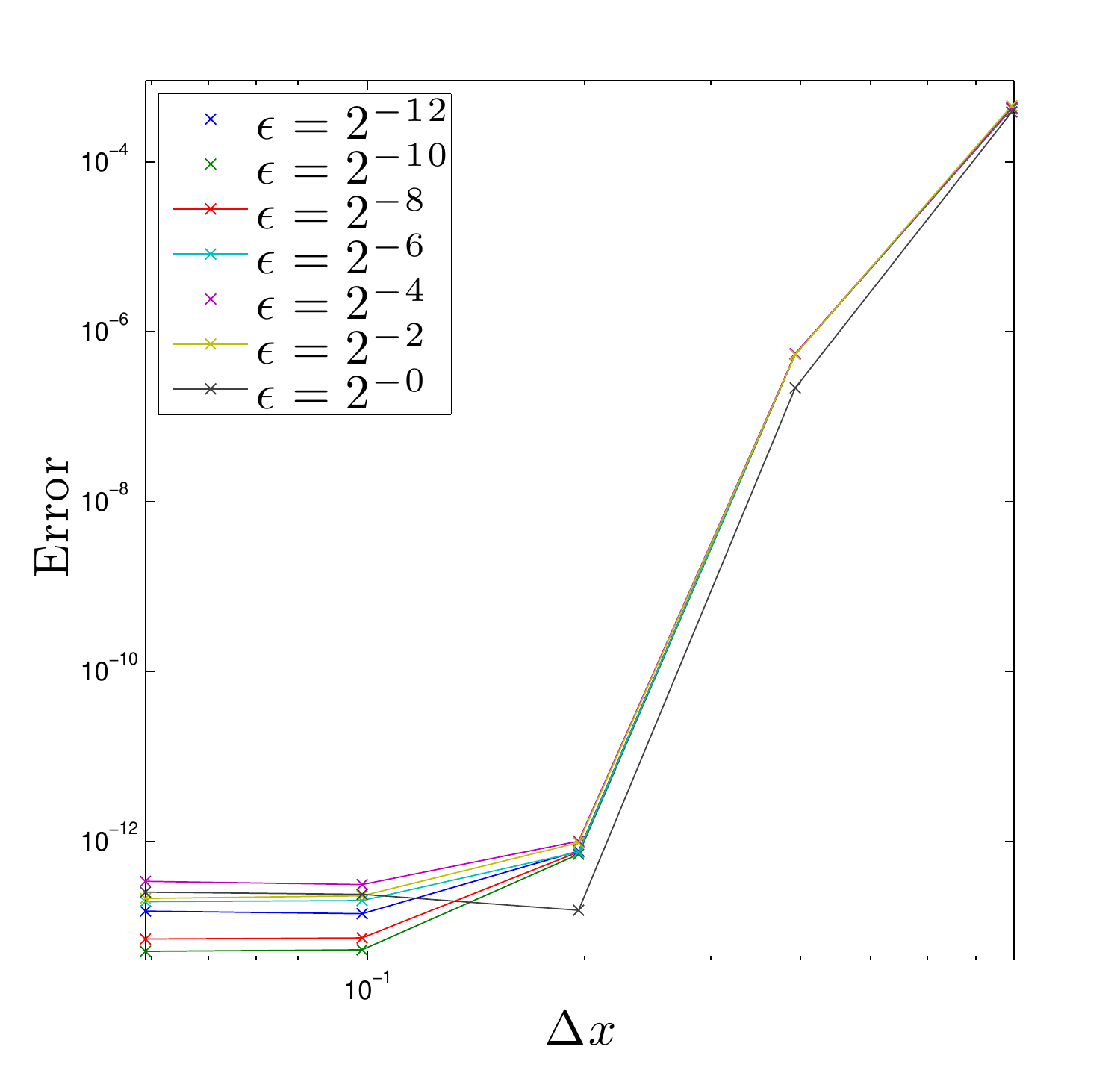}
		\caption{\label{fig:Tf1cvgdxord4v1} $err_{\rho^\eps} (T_f = 0.1)$ w.r.t $\Delta x$, $N_t = 2^{10}$} 
		\end{subfigure}
		\hspace{-1cm}
		\begin{subfigure}[t]{0.6\textwidth}
		\centering
		\includegraphics[height=5.cm,width=\textwidth]{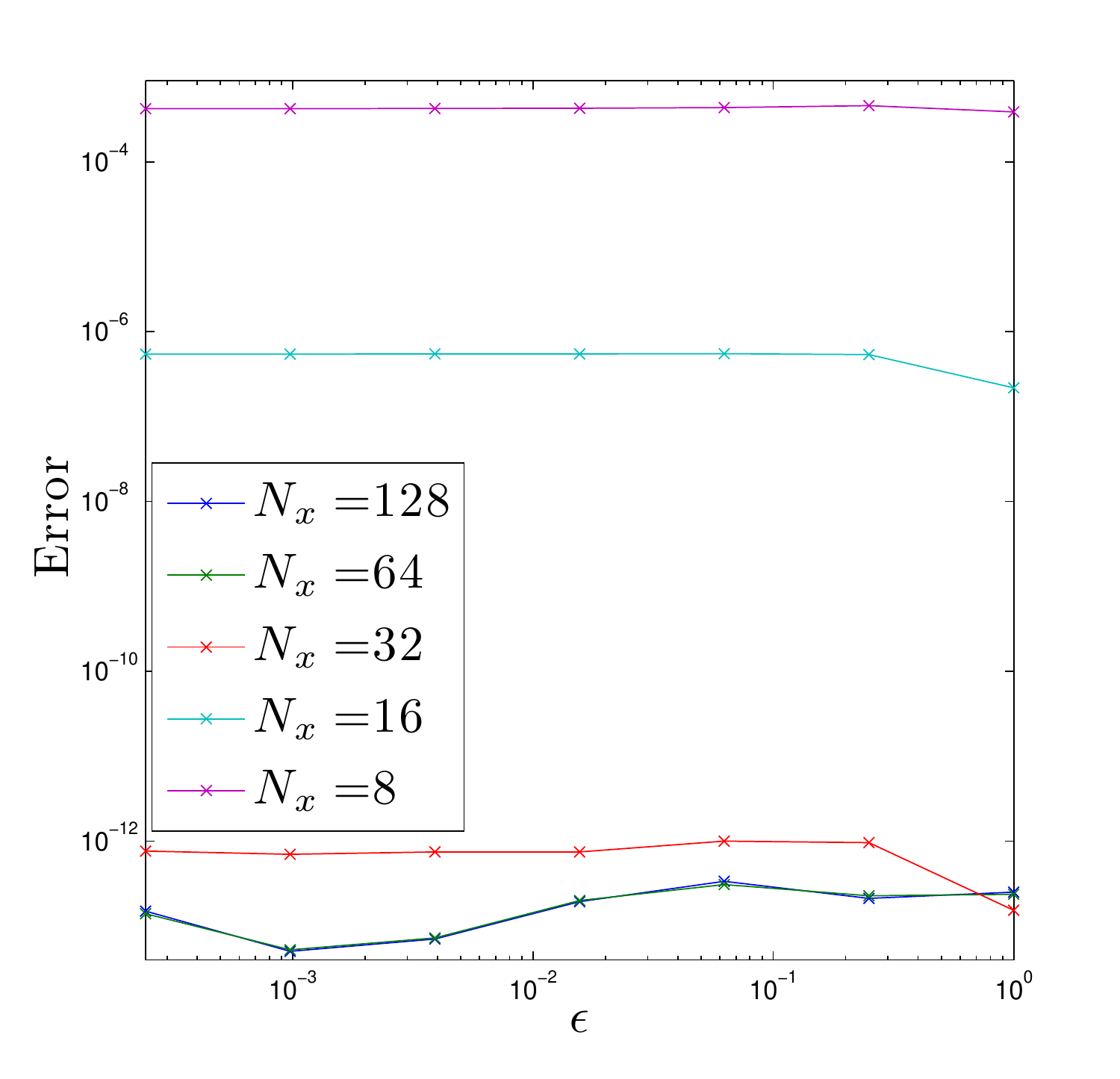}
		\caption{\label{fig:Tf1cvgdxord4v2} $err_{\rho^\eps} (T_f = 0.1)$ w.r.t $\eps$, $N_t = 2^{10}$} 
		\end{subfigure}
		\hspace{-6cm}
		\caption{\label{fig:Tf1cvgdxord4-1} Error on $\rho^\eps$ for the splitting scheme \eqref{scheme4} of order $4$ before the caustics: dependence on $\eps$ and on $\Delta x$.}
		\end{figure}

\begin{figure}[p]
		\centering	
		\hspace{-5cm}
		\begin{subfigure}[t]{0.6\textwidth}
		\centering
		\includegraphics[height=5.cm,width=\textwidth]{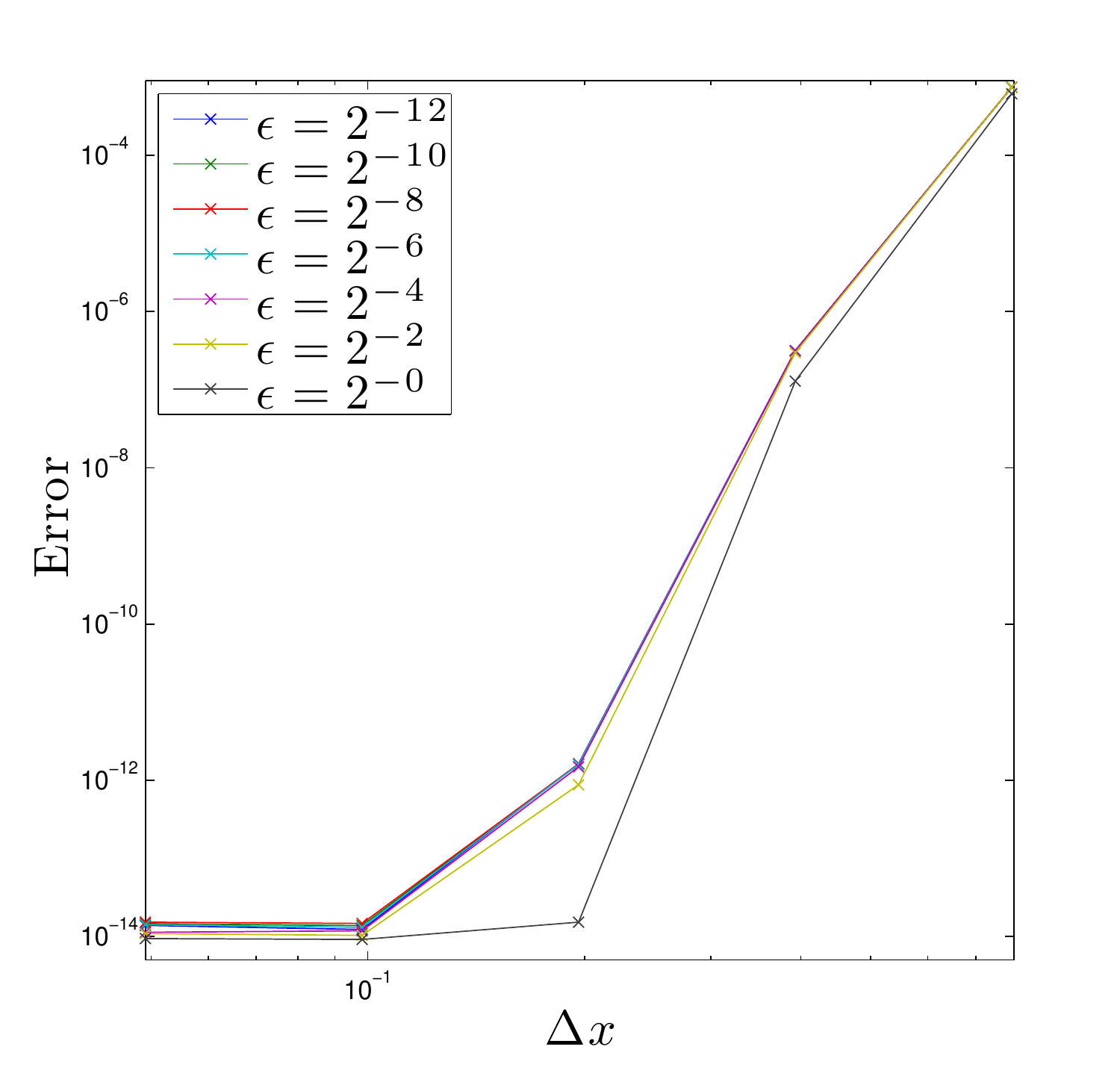}
		\caption{\label{fig:Tf1cvgdxord4v3} $err_{(\se,\ae)} (T_f = 0.1)$ w.r.t $\Delta x$, $N_t = 2^{10}$} 
		\end{subfigure}
		\hspace{-1cm}
		\begin{subfigure}[t]{0.6\textwidth}
		\centering
		\includegraphics[height=5.cm,width=\textwidth]{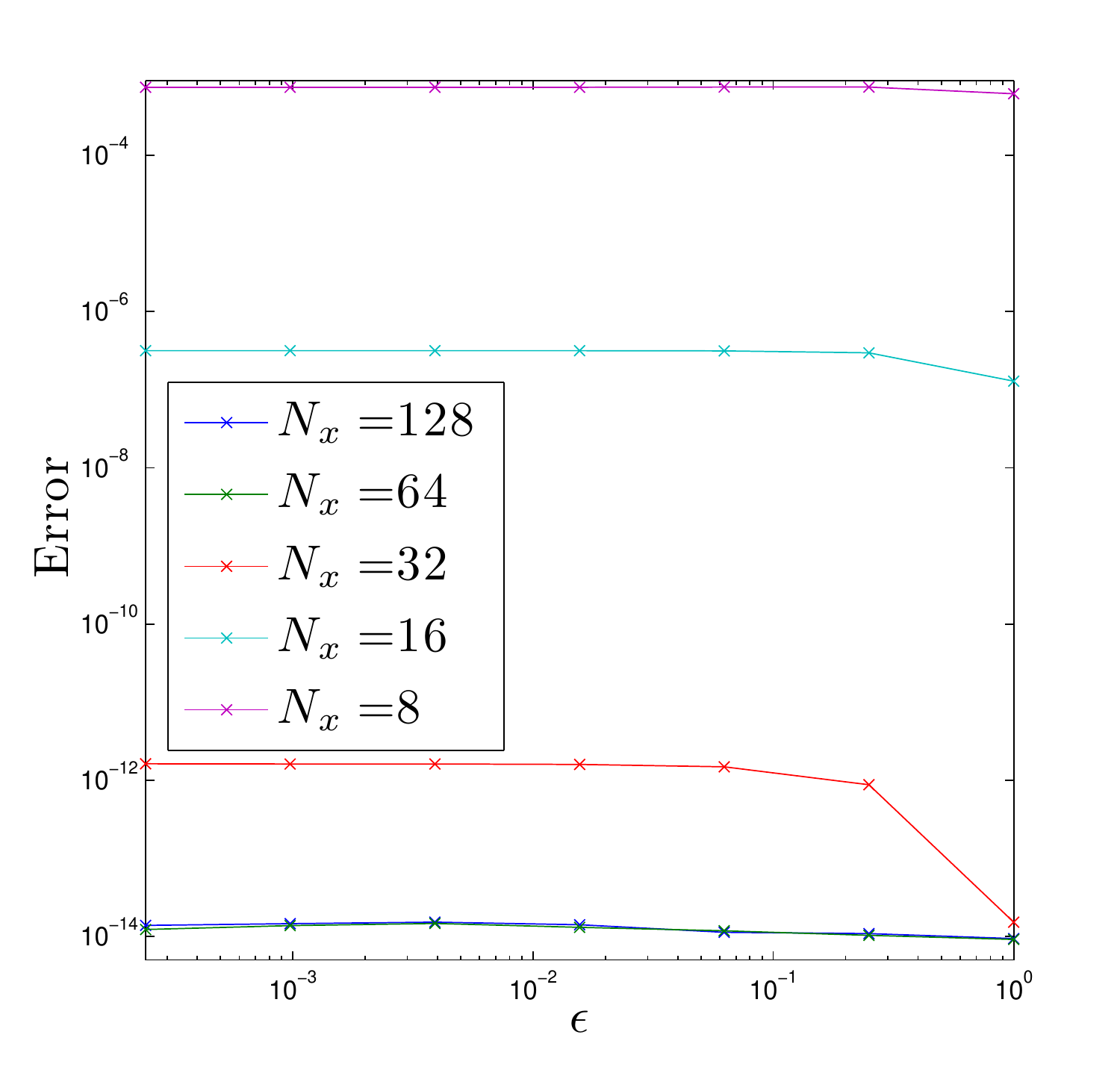}
		\caption{\label{fig:Tf1cvgdxord4v4} $err_{(\se,\ae)} (T_f = 0.1)$ w.r.t $\eps$, $N_t = 2^{10}$} 
		\end{subfigure}
		\hspace{-6cm}
				\caption{\label{fig:Tf1cvgdxord4-2}Error on $(\se,\ae)$ for the splitting scheme \eqref{scheme4} of order $4$ before the caustics: dependence on $\eps$ and on $\Delta x$.}

		\end{figure}

	Figures \ref{fig:Tf1cvgtpsord4-1} to \ref{fig:Tf1cvgdxord4-2} illustrate the behavior of our fourth-order scheme \eqref{scheme4} at $T_f = 0.1$: here again, it appears that, before the caustics, our method is UA with an order 4 in time and with spectral in space accuracy.
	
\begin{figure}[p]
	\centering	
		\hspace{-5cm}
		\begin{subfigure}[t]{0.7\textwidth}
		\centering
		\includegraphics[height=5.cm,width=\textwidth]{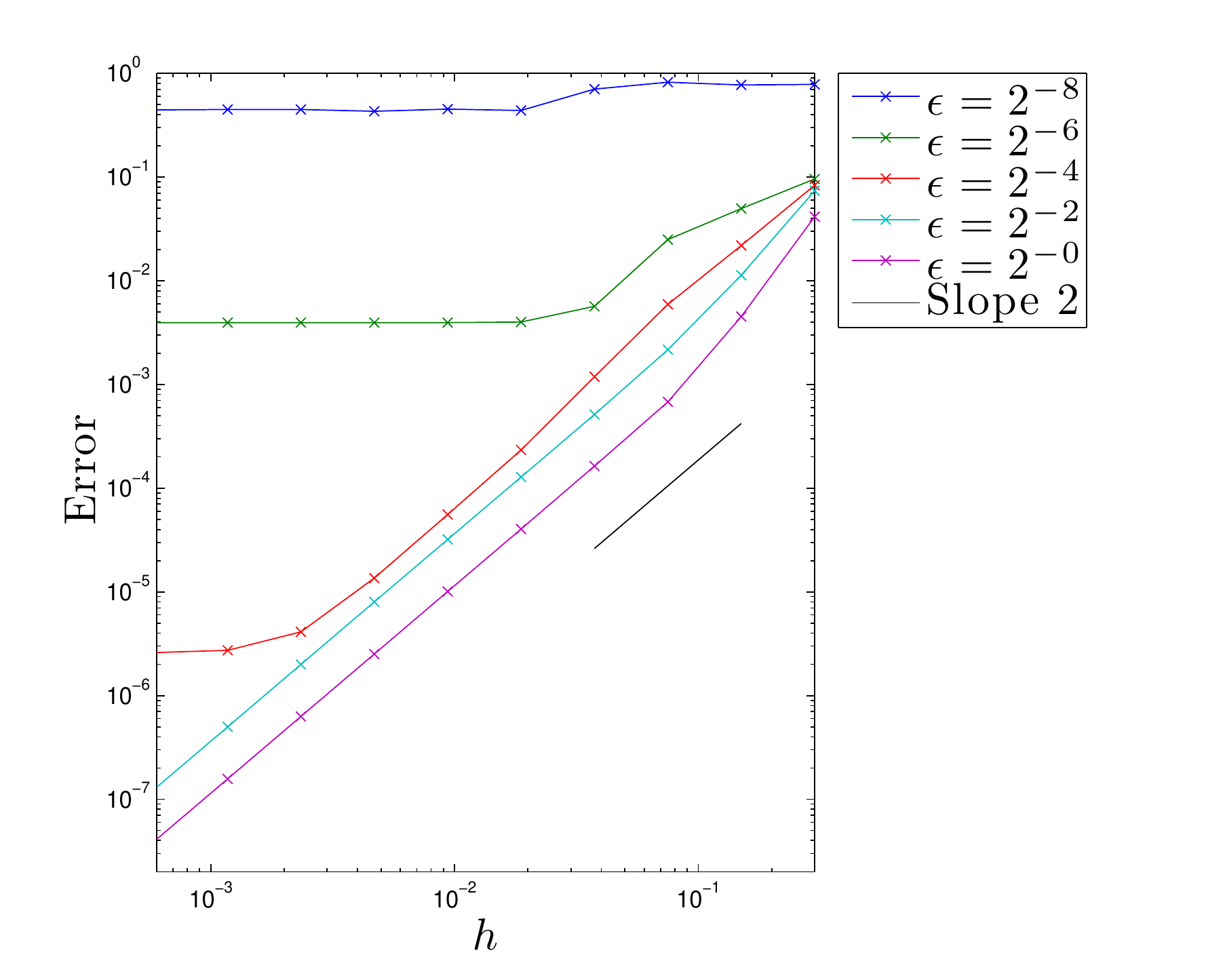}
		\caption{\label{fig:Tf6cvgtpsord2v1NLS} $err_{\rho^\eps} (T_f = 0.6)$ w.r.t $h$, $N_x = 2^7$ } 
		\end{subfigure}
		\hspace{-1cm}
		\begin{subfigure}[t]{0.7\textwidth}
		\centering
		\includegraphics[height=5.cm,width=\textwidth]{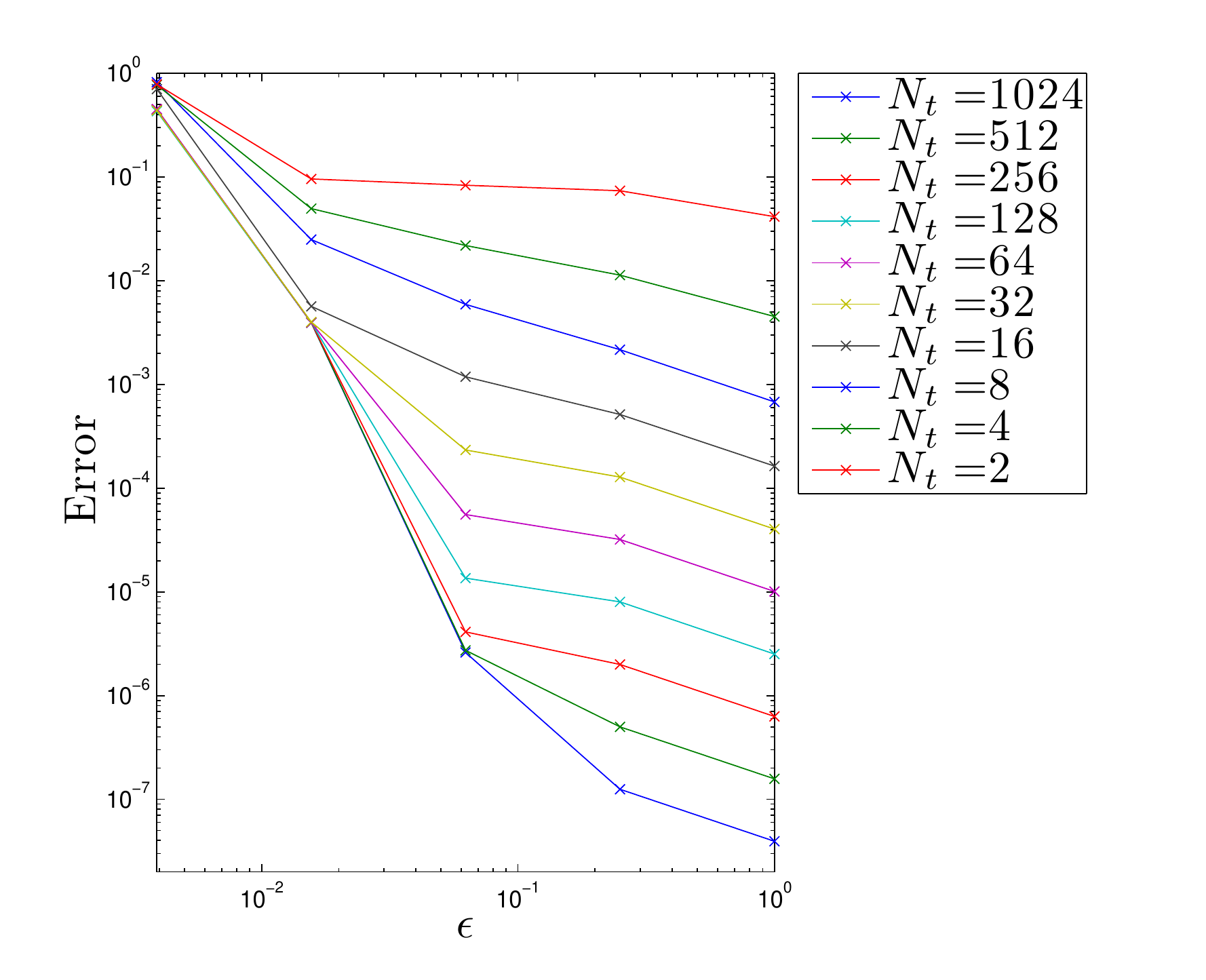}
		\caption{\label{fig:Tf6cvgtpsord2v2NLS} $err_{\rho^\eps} (T_f = 0.6)$ w.r.t $\eps$, $N_x = 2^7$} 
		\end{subfigure}
		\hspace{-6cm}
		\caption{ \label{fig:Tf6cvgtpsord2NLS}Error on $\rho^\eps$ for the Strang splitting scheme for \eqref{eq:GPE} after the caustics, dependence on $\eps$ and on $h$.}
	
		\end{figure}

\begin{figure}[p]
		\centering	
		\hspace{-5cm}
		\begin{subfigure}[t]{0.7\textwidth}
		\centering
		\includegraphics[height=5.cm,width=\textwidth]{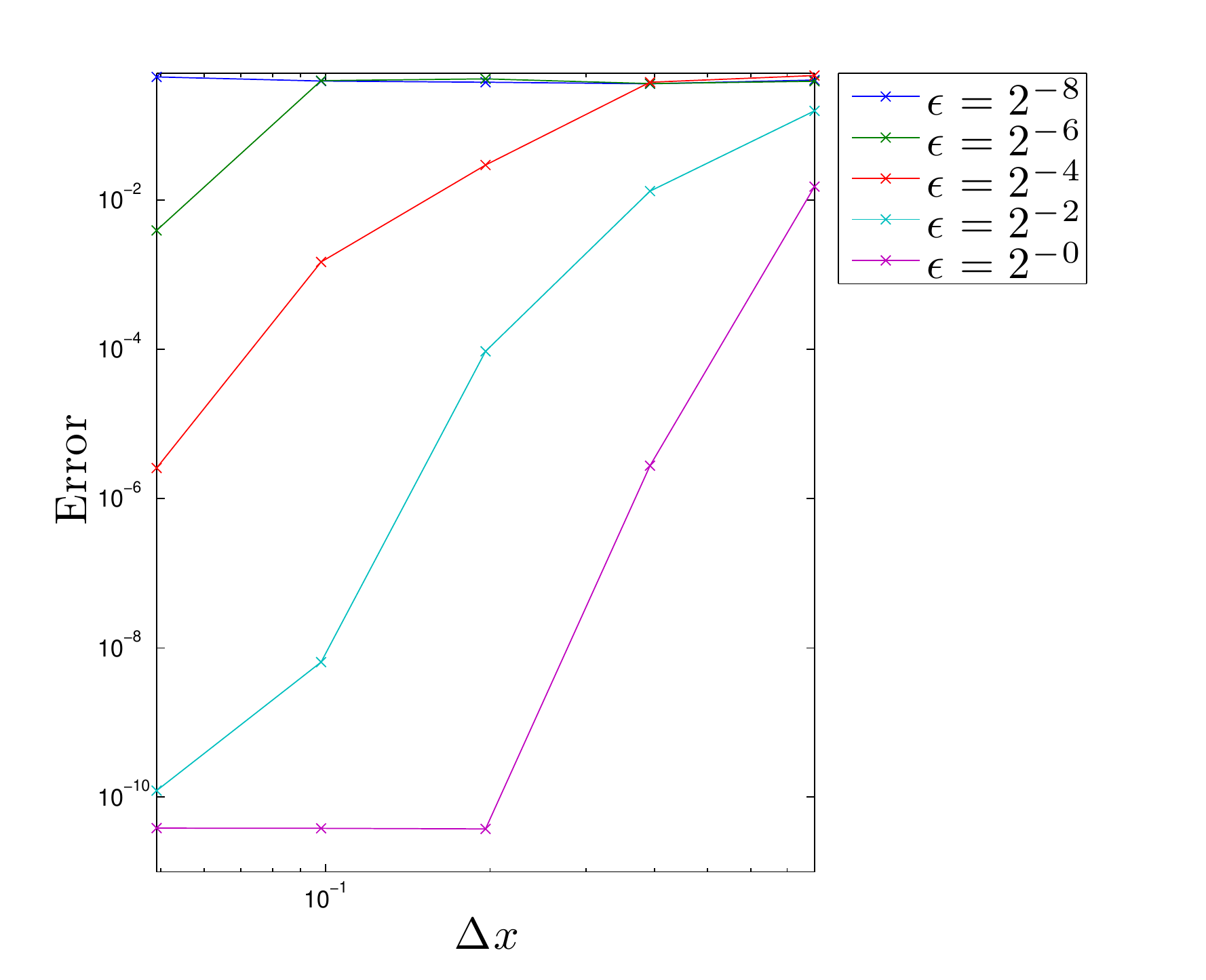}
		\caption{\label{fig:Tf6cvgdxord2v1NLS} $err_{\rho^\eps} (T_f = 0.6)$ w.r.t $\Delta x$, $N_t = 2^{15}$  } 
		\end{subfigure}
		\hspace{-1cm}
		\begin{subfigure}[t]{0.7\textwidth}
		\centering
		\includegraphics[height=5.cm,width=\textwidth]{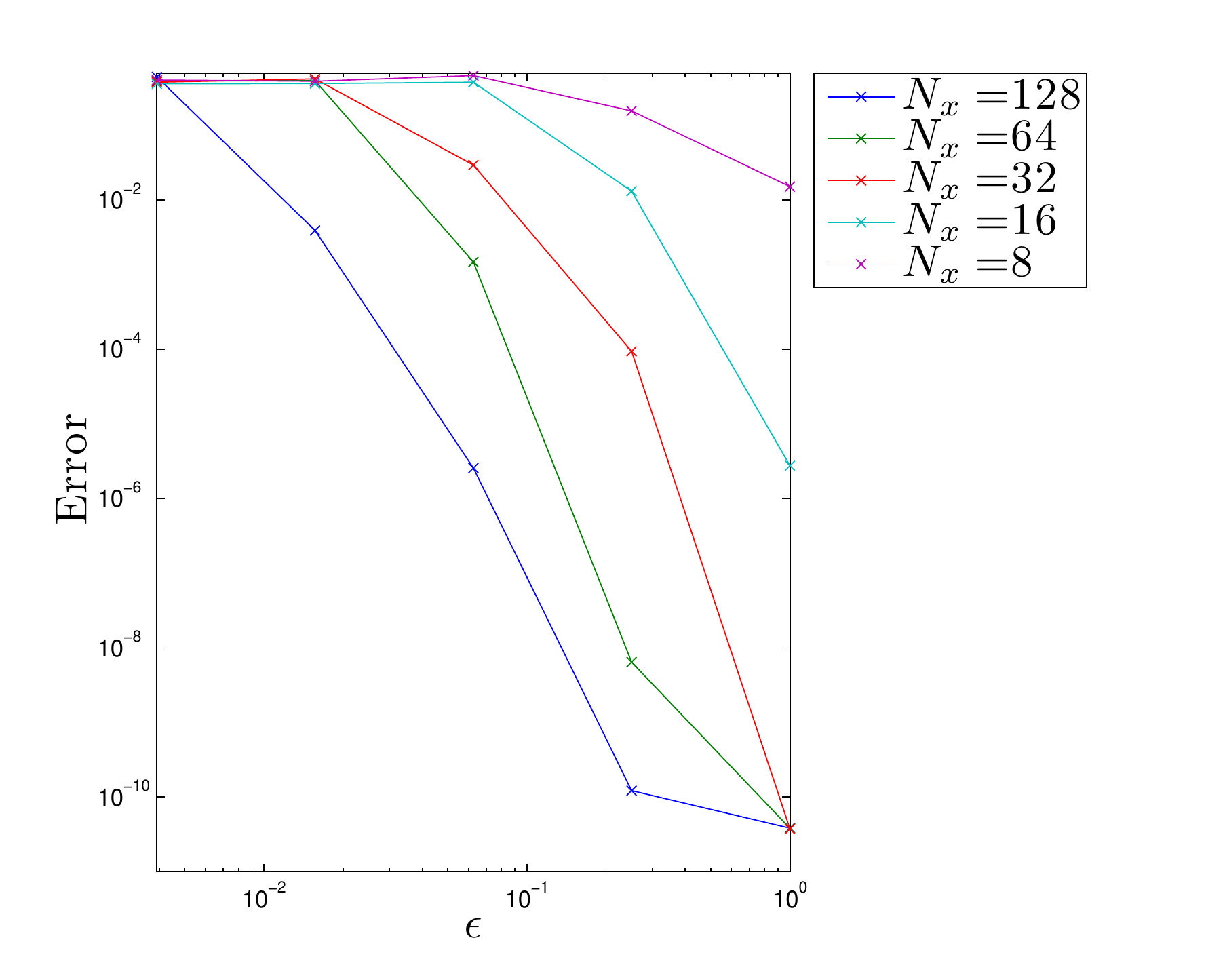}
		\caption{\label{fig:Tf6cvgdxord2v2NLS} $err_{\rho^\eps} (T_f = 0.6)$ w.r.t $\Delta x$, $N_t = 2^{15}$ } 
		\end{subfigure}
		\hspace{-6cm}
		\caption{\label{fig:Tf6cvgdxord2NLS}Error on $\rho^\eps$ for the Strang splitting scheme for \eqref{eq:GPE} after the caustics, dependence on $\eps$ and on $\Delta x$.}

		\end{figure}

\begin{figure}[p]
	\centering	
		\hspace{-5cm}
		\begin{subfigure}[t]{0.6\textwidth}
		\centering
		\includegraphics[height=5.cm,width=\textwidth]{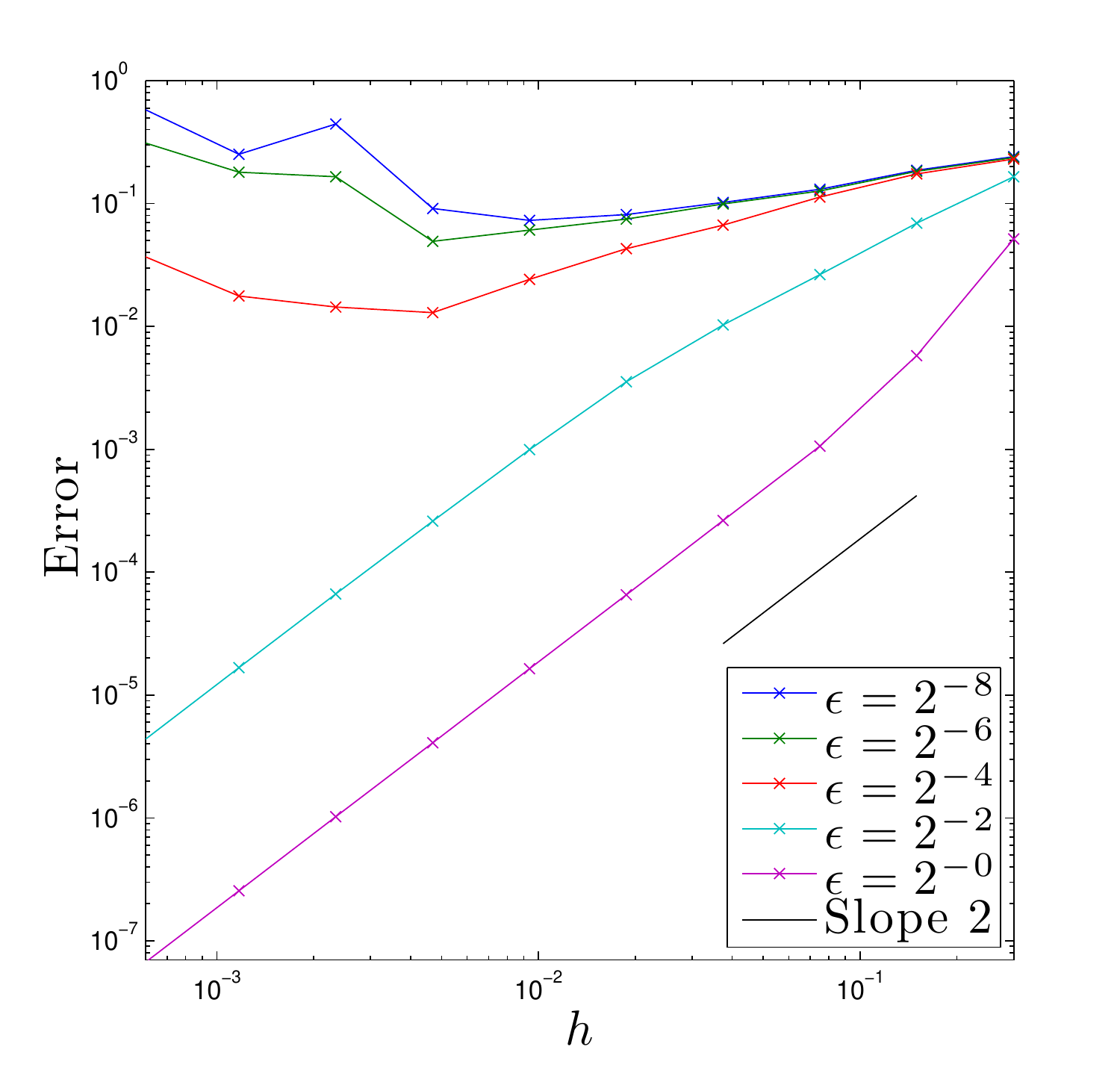}
		\caption{\label{fig:Tf6cvgtpsord4v1} $err_{\rho^\eps} (T_f = 0.6)$ w.r.t $h$, $N_x = 2^7$ } 
		\end{subfigure}
		\hspace{-1cm}
		\begin{subfigure}[t]{0.7\textwidth}
		\centering
		\includegraphics[height=5.cm,width=\textwidth]{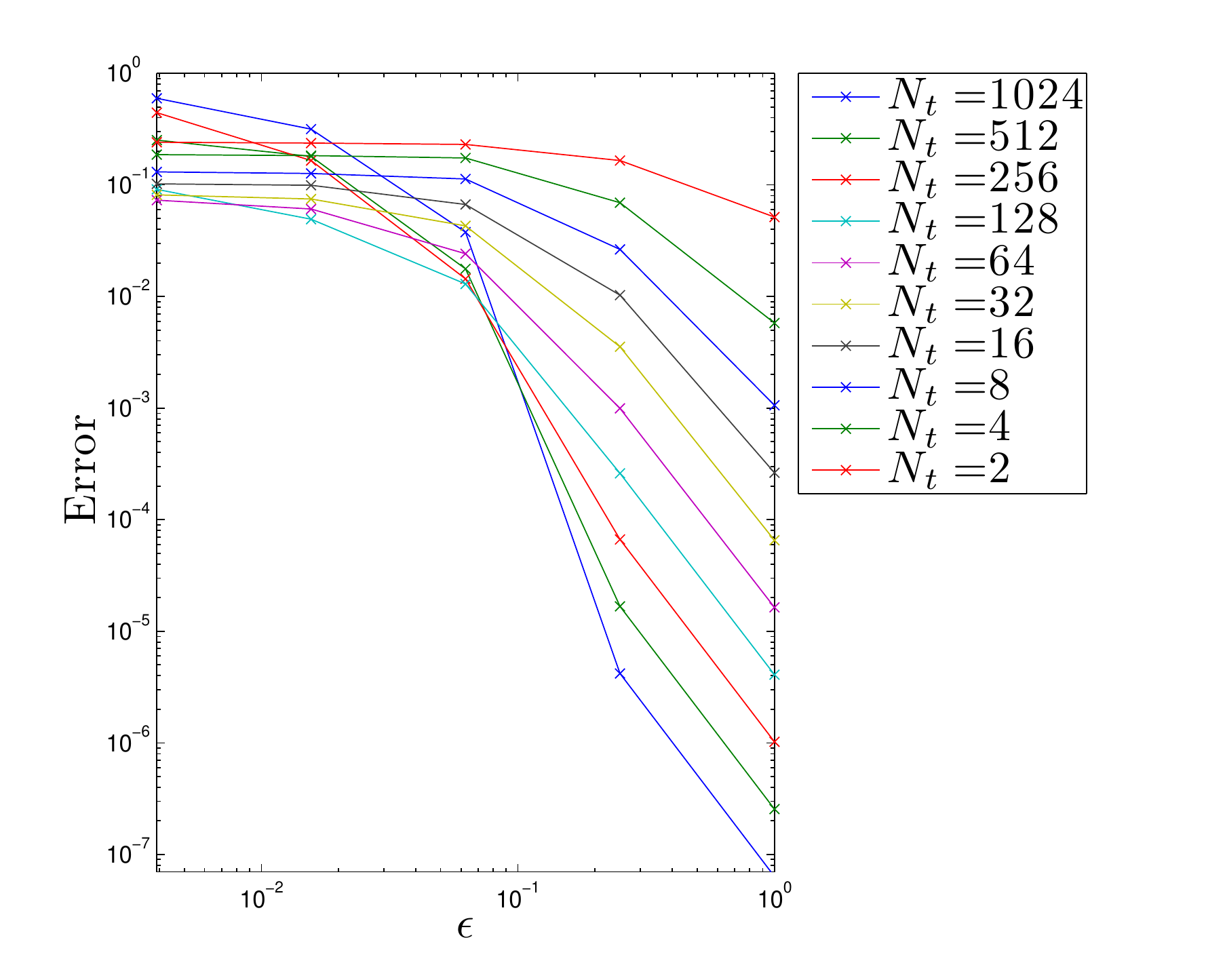}
		\caption{\label{fig:Tf6cvgtpsord4v2} $err_{\rho^\eps} (T_f = 0.6)$ w.r.t $\eps$, $N_x = 2^7$} 
		\end{subfigure}
		\hspace{-6cm}
		\caption{ \label{fig:Tf6cvgtpsord4}Error on $\rho^\eps$ for the splitting scheme \eqref{scheme2} of order $2$ after the caustics, dependence on $\eps$ and on $h$.}
	
		\end{figure}

\begin{figure}[p]
		\centering	
		\hspace{-5cm}
		\begin{subfigure}[t]{0.7\textwidth}
		\centering
		\includegraphics[height=5.cm,width=\textwidth]{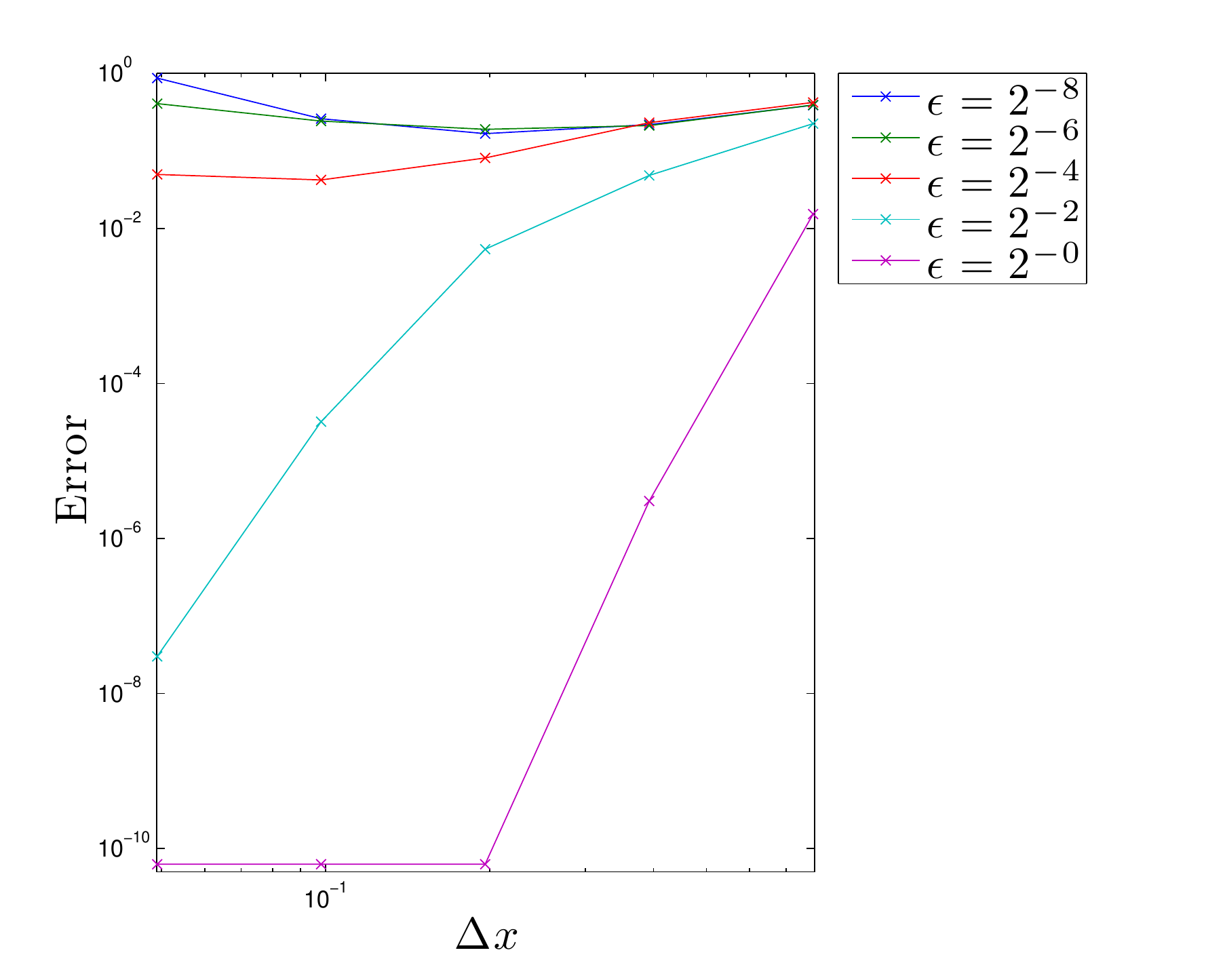}
		\caption{\label{fig:Tf6cvgdxord2v1} $err_{\rho^\eps} (T_f = 0.6)$ w.r.t $\Delta x$, $N_t = 2^{15}$  } 
		\end{subfigure}
		\hspace{-1cm}
		\begin{subfigure}[t]{0.6\textwidth}
		\centering
		\includegraphics[height=5.cm,width=\textwidth]{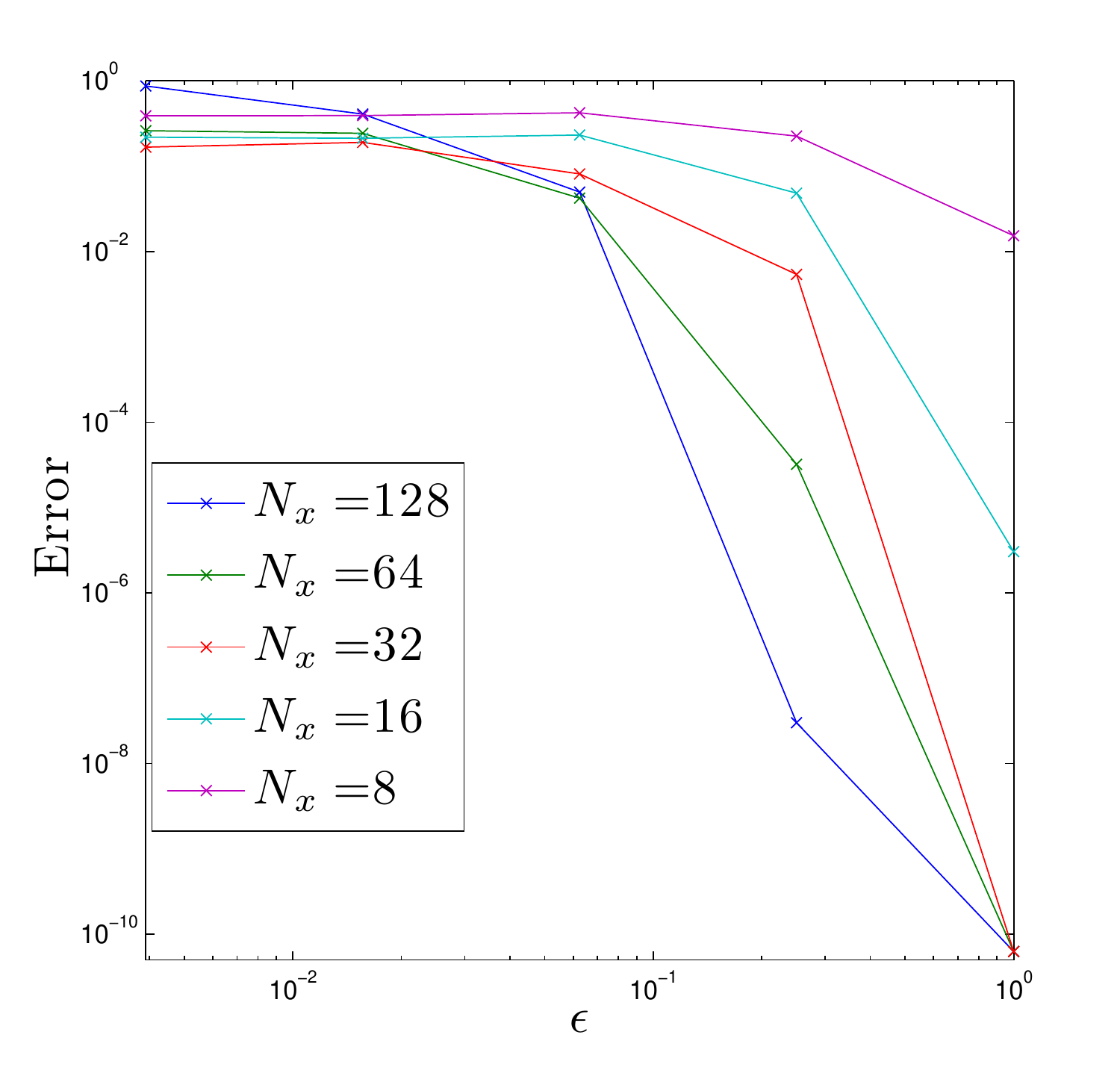}
		\caption{\label{fig:Tf6cvgdxord2v2} $err_{\rho^\eps} (T_f = 0.6)$ w.r.t $\eps$, $N_t = 2^{15}$ } 
		\end{subfigure}
		\hspace{-6cm}
		\caption{\label{fig:Tf6cvgdxord4}Error on $\rho^\eps$ for the splitting scheme \eqref{scheme2} of order $2$ after the caustics, dependence on $\eps$ and on $\Delta x$.}

		\end{figure}

	Finally, let us explore the behavior of the splitting methods after caustics, by observing the error on the density $\rho^\eps$. Figures \ref{fig:Tf6cvgtpsord2NLS}, \ref{fig:Tf6cvgdxord2NLS}, \ref{fig:Tf6cvgtpsord4} and \ref{fig:Tf6cvgdxord4} present the same simulations as Figures \ref{fig:Tf1cvgtpsord2-1NLS}, \ref{fig:Tf1cvgdxord2-1NLS}, \ref{fig:Tf1cvgtpsord2-1} and \ref{fig:Tf1cvgdxord2-1}, except that the final time is now $T_f = 0.6$, \textit{i.e.} we illustrate the behaviors of Strang splitting method and of scheme \eqref{scheme2} after the caustics. In that case, it appears that none of these methods is UA, neither in $h$, nor in $\Delta x$, with respect to $\eps$. Concerning the Strang splitting scheme, this behavior was already reported in \cite{MR2359915,MR2805153}. Notice that, although it is not UA any longer, our scheme \eqref{scheme2} still has second-order accuracy in time and spectral accuracy in space (with $\eps$-dependent constants). Recall that the same scheme written on \eqref{eq:SCv1} would not be usable in the same situation, since $S^\eps$ ceases to be regular for $\eps>0$, after the formation of caustics.
	
	All the experiments have been performed with the methods of both Sections \ref{sec:semilag} and \ref{sec:splieik} for the eikonal equation  \eqref{eq:phase2} but we just present here the graphs obtained with the semi-lagrangian methods of Section \ref{sec:semilag}.
	
\section{Final remarks}
Let us emphasize that, as a by-product, we have also derived a new numerical scheme based on splitting techniques to approximate the solution of the Hamilton-Jacobi (eikonal) equation
\be\label{eq:phase2}
	\pa_t S + \frac{|\nabla S|^2}{2} = 0
\ee
based on the Cole-Hopf transform. Finally, it is interesting to mention that, although  we have chosen to focus here on the supercritical regime of \eqref{eq:GPE} (see \cite{bookRemi}), our approach is also relevant in other semiclassical regimes, whether the linear Schr\"odinger equation with a given potential, or the weakly nonlinear geometric optics, where  $|\pe|^2\pe$ is replaced with $\eps|\pe|^2\pe$.

		

\bibliographystyle{siam}
\bibliography{bibliographiebibdesk}		
 \end{document}